\documentclass{article}
\newcommand{\upr}{\overline{P}}
\newcommand{\lpr}{\underline{P}}

\newcommand{\lpq}{\underline{Q}}
\newcommand{\upe}{\overline{E}}
\newcommand{\lpe}{\underline{E}}

\newcommand{\lreg}{\underline{R}}

\newcommand{\rset}{\mathbb{R}}

\newcommand{\impev}{\emptyset}

\newcommand{\lset}{\mathcal{L}}

\newcommand{\prt}{I\dsn P}

\newcommand{\dset}{\mathcal{D}}
\newcommand{\aset}{\mathcal{A}}

\newcommand{\mset}{\mathcal{M}}

\newcommand{\bset}{\mathcal{B}}

\newcommand{\asetpa}{\aset(\prt)}

\newcommand{\dsn}{\!\!}

\newcommand{\comment}[1]{}
\newcommand{\nega}[1]{{#1}^\mathsf{c}}
{\left\lbrace\begin{array}{@{}l@{}}}%
	{\end{array}\right.}

\usepackage{amssymb}
\usepackage{amsmath}
\usepackage{amsthm}
\usepackage[nice]{nicefrac}
\usepackage{xfrac}
\usepackage{mathtools}
\usepackage{authblk}
\usepackage{graphicx}
\usepackage[misc]{ifsym}
\usepackage{bbding}
\usepackage{lineno}
\newtheorem{theorem}{Theorem}[section]
\newtheorem{proposition}{Proposition}[section]
\newtheorem{lemma}{Lemma}[section]
\newtheorem{corollary}{Corollary}[section]
\newtheorem{remark}{Remark}[section]
\newtheorem{definition}{Definition}[section]
\newtheorem{example}{Example}[section]

\begin{document}
	
\title{Dilation Properties of Coherent Nearly-Linear Models}

\author[1]{Renato Pelessoni\thanks{renato.pelessoni@deams.units.it}}
\author[1]{Paolo Vicig\thanks{paolo.vicig@deams.units.it}}
\affil[1]{DEAMS ``B. de Finetti''\\
	University of Trieste\\
	Piazzale Europa~1\\
	I-34127 Trieste\\
	Italy}

\renewcommand\Authands{ and }

\maketitle

\begin{abstract}
Dilation is a puzzling phenomenon within Imprecise Probability theory: when it obtains, our uncertainty evaluation on event $A$ is vaguer after conditioning $A$ on $B$, \emph{whatever} is event $B$ in a given partition $\bset$.
In this paper we investigate dilation with coherent Nearly-Linear (NL) models.
These are a family of neighbourhood models, obtaining lower/upper probabilities by linear affine transformations (with barriers) of a given probability, and encompass several well-known models, such as the Pari-Mutuel Model, the $\varepsilon$-contamination model, the Total Variation Model, and others.
We first recall results we recently obtained for conditioning NL model with the standard procedure of natural extension and separately discuss the role of the alternative regular extension.
Then, we characterise dilation for coherent NL models.
For their most relevant subfamily, Vertical Barrier Models (VBM), we study the coarsening property of dilation, the extent of dilation, and constriction.
The results generalise existing ones established for special VBMs.
As an interesting aside, we discuss in a general framework how logical (in)dependence of $A$ from $\bset$ or extreme evaluations for $A$ influence dilation.	


\end{abstract}

\section*{Acknowledgement}
*NOTICE: This is the authors' version of a work that was accepted for publication in International Journal of Approximate Reasoning. Changes resulting from the publishing process, such as peer review, editing, corrections, structural formatting, and other quality control mechanisms may not be reflected in this document. Changes may have been made to this work since it was submitted for publication. A definitive version was subsequently published in International Journal of Approximate Reasoning, 
vol. 140, January~2022, pages 211–-231
https://doi.org/10.1016/j.ijar.2021.10.009 $\copyright$ Copyright Elsevier

\vspace{0.3cm}
$\copyright$ 2019. This manuscript version is made available under the CC-BY-NC-ND 4.0 license http://creativecommons.org/licenses/by-nc-nd/4.0/

\begin{center}
	\includegraphics[width=2cm]{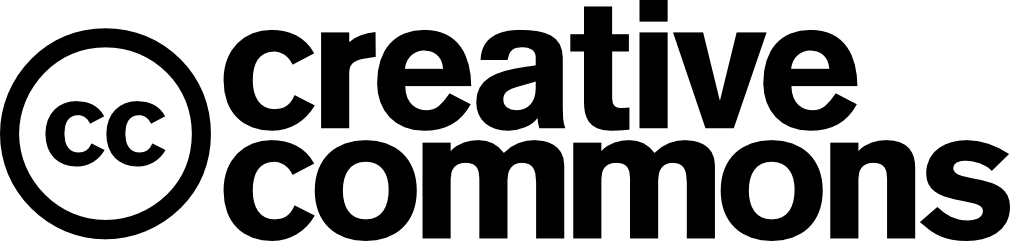}
	\includegraphics[width=2cm]{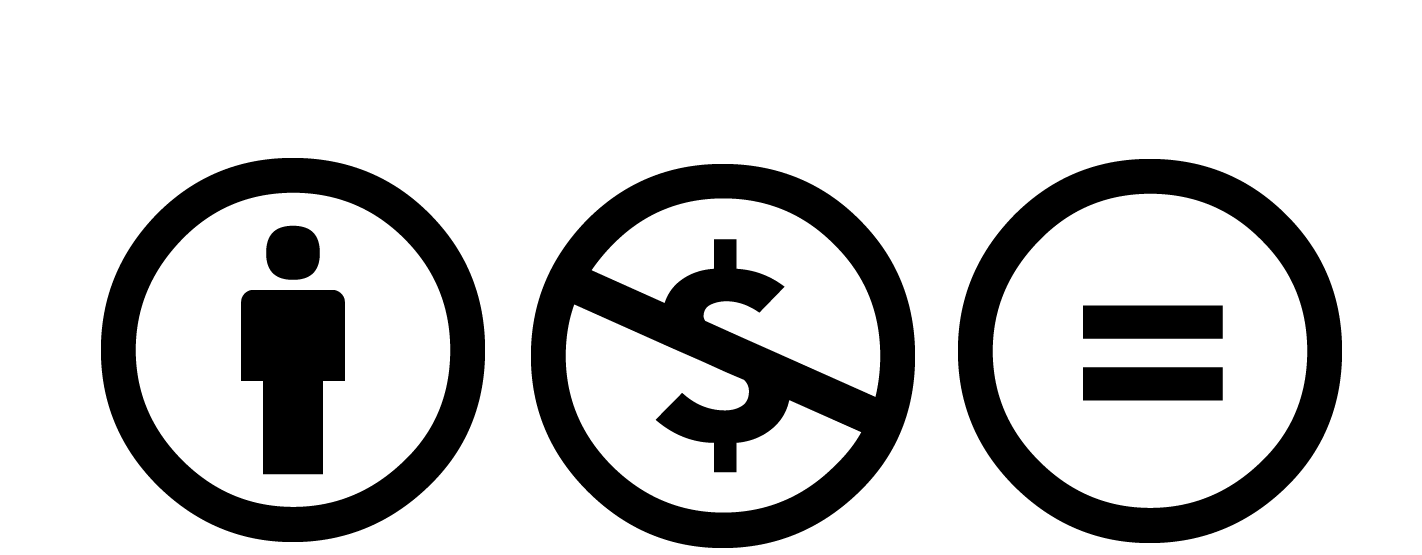}
\end{center}

\section{Introduction}
\label{sec:introduction}
Several uncertainty measurement models are currently grouped under the name \emph{Imprecise Probabilities}, ranging from very general ones, like coherent lower and upper previsions, to others more specific, like possibilities and necessities.
A special class of models, relatively simple to work with, is that of \emph{neighbourhood models} \cite{walley_statistical_1991} or distortion models \cite{montes_neighbourhood_I, montes_neighbourhood_II}.
Here a lower ($\lpr$) and an upper ($\upr$) probability are obtained as functions of a given (precise) probability $P_0$.
Well known examples are the Pari-Mutuel Model (PMM) \cite{montes_pari-mutuel_2019,pelessoni_inference_2010,walley_statistical_1991}, the $\varepsilon$-contamination model or linear-vacuous mixture \cite{herron_divisive_1997,walley_statistical_1991}, the Total Variation Model (TVM) \cite{herron_divisive_1997} and others.
In previous work \cite{corsato_nearly-linear_2019} we introduced a family of such neighbourhood models,
termed \emph{Nearly-Linear} (NL) models, where $\lpr$ and $\upr$ are linear affine transformations of $P_0$, with barriers to guarantee that $\lpr$, $\upr$ belong to the $[0,1]$ interval.
Despite their functional simplicity, NL models can formalise a certain variety of beliefs, including some (common, but) weakly consistent ones.
In fact, they are partitioned into three subfamilies, of which that of \emph{Vertical Barrier Models} (VBM) is formed by $(\lpr,\upr)$ that are coherent, while in the other two subfamilies only Horizontal Barrier Models (HBM) may be coherent in non-trivial cases, but they generally satisfy weaker consistency requirements.
By contrast, VBMs include the known coherent models recalled above (PMM, $\varepsilon$-contamination, TVM) and others as subcases.

We investigated several features of NL models, as well as their relationships to other models, in \cite{corsato_nearly-linear_2019, pelessoni_dilation_2020, pelessoni_inference_2020}.
In particular, in \cite{pelessoni_dilation_2020} we obtained early results concerning conditioning and dilation with coherent NL models, which are developed and extended in this paper.

Precisely, we first recall some preliminary notions in Section \ref{sec:preliminaries}.
These include some minimal information on \emph{Williams' coherence} \cite{williams_notes_2007}, which is the consistency notion we apply in a conditional environment, and the basic inferential procedure termed \emph{natural extension} \cite{pelessoni_williams_2009,troffaes_lower_2014,walley_statistical_1991}. 
Section \ref{subsub:Nearly_Linear_Models} describes briefly NL models, whose $\lpr$, $\upr$ are defined on the set $\asetpa$ of events logically dependent on a partition $\prt$ (i.e. on the powerset of $\prt$).

Conditioning a coherent NL model on a given event $B\in\asetpa\setminus\{\emptyset\}$ is discussed in Section \ref{sec:NLM_conditioning}.
We apply conditioning with the natural extension, already investigated in \cite{pelessoni_dilation_2020}. It makes use of simple formulae available for $2$-monotone ($2$-alternating) lower (upper) probabilities, such as the imprecise probabilities making up coherent NL models.
These formulae further specialise with VMBs, showing that the resulting model is still a VBM, a property (stability) not shared by HBMs.

In the next sections we focus on investigating \emph{dilation} with coherent NL models.
Consider conditioning some $A\in\asetpa$ on any event $B$ in a partition $\bset$.
Dilation occurs when both $\lpr(A|B)\leq\lpr(A)$ and $\upr(A|B)\geq\upr(A)$, for any $B\in\bset$.
Dilation is known to be a phenomenon not really uncommon with imprecise probabilities \cite{walley_statistical_1991}, and has been studied in dedicated papers, among which \cite{herron_divisive_1997,seidenfeld_dilation_1993} and more recently \cite{nielsen_counterexamples_2019,pedersen_dilation_2015}.
For the present paper, \cite{herron_divisive_1997,seidenfeld_dilation_1993} are particularly significant, because they, making use (implicitly) of the natural extension, $(a)$ investigate dilation for some special VBMs, $(b)$ define properties of dilation such as the coarsening property or the extent of dilation.

In Section \ref{sec:dilation_NL}, after recalling definitions involving dilation (Section \ref{subsec:what_dilation}), the results from Section \ref{sec:dilation_generalities} in the Appendix contribute to fixing the basic assumptions for investigating dilation with NL models in Section \ref{subsec:dilation_with_NL}.
With these assumptions, our framework remains more general than in \cite{herron_divisive_1997,seidenfeld_dilation_1993}.
In particular, it is possible that $\lpr(B)=0$ for $B\in\bset$.

We discuss the role of logical independence of $A$ from $\bset$ with NL models in Section \ref{sec:role_log_ind}.
Not surprisingly, logical independence limits the number of events $A\in\asetpa$ that $\bset$ may dilate, and tends to vanish by increasing the cardinality of partition $\bset$.
In Section \ref{sec:results_dilation} we give a characterisation of dilation for coherent NL models (Proposition \ref{pro:dilation_NLM}) and show (Proposition \ref{dilation-epsilon}) that for the $\varepsilon$-contamination model it boils down to the characterisation in \cite{herron_divisive_1997}.

In Section \ref{sec:properties_dilation_VBM} we discuss properties of dilation with reference to VBMs.
Section \ref{sec:coarsening} concerns the \emph{coarsening property}, that obtains if knowing that $\bset$ dilates $A$ implies  that there exists a partition coarser than $\bset$ that dilates $A$.
The main result is Theorem \ref{thm:coarsening_basic}, a sufficient condition for coarsening, extending achievements in \cite{herron_divisive_1997}.
In Section \ref{subsec:extent_dilation} the \emph{extent of dilation} is computed (Theorem \ref{thm:extent_dilation_VBM}).
A sufficient condition for imprecision increase is derived in Proposition \ref{pro:impr_incr_bmin1}.
Section \ref{subsec:constriction} discusses \emph{constriction}, a sort of opposite of dilation.
Constriction is a very desirable property, but unfortunately very hard to obtain, as appears from the results in this section.
Finally, our conclusions are presented in Section \ref{sec:conclusions}.

Two questions, relevant in general to the investigation of dilation, are discussed in the Appendix.
Its Section~\ref{subsec:cond_RE} is concerned with the first question, the role of the \emph{regular extension}.
This is an alternative conditioning aiming at limiting
those instances where $\lpr(B)=0$ induces vague inferences via natural extension.
The regular extension was introduced in \cite{walley_coherent_1981} and discussed in \cite[Appendix J]{walley_statistical_1991} and other papers, including \cite{miranda_coherent_2015}, referring to Walley's coherence notions,
while results concerning its Williams' coherence are only indirectly available.
In the Appendix, we approach the regular extension 
$(a)$ proving directly its Williams' coherence in a general framework (Section \ref{subsub:W-coh_reg_ext}), $(b)$ showing that within $2$-monotone ($2$-alternating) models, such as coherent NL models in particular, it differs from the natural extension under quite restrictive conditions (Section \ref{subsub:reg_ext_2_mon_mod}), and $(c)$ that for VBMs these conditions imply that the VBM is a Pari-Mutuel Model (Section \ref{subsub:case_VBM}).
This lets us conclude that the regular extension plays a rather limited role with coherent NL models.

Section~\ref{sec:dilation_generalities}, the second part of the Appendix, tackles the other question, not 
explicitly discussed in the literature, of which assumptions to require in investigating dilation.
Clearly, hypotheses that let us know \emph{a priori} whether it occurs or not, and only these, should be discarded.
We analyse this problem in general (not restricting to NL models only),
showing $(a)$ that logical independence of event $A$ from partition $\bset$ should obtain (Section \ref{subsubsec:dil_log_ind}) and $(b)$ that, while conditioning with the natural extension, $A$ should not have precise probability $0$ or $1$ (Section~\ref{subsubsec:dilation_nat_ext_extr_ev}).

\section{Preliminaries}
\label{sec:preliminaries}
In this paper we shall be concerned with coherent lower and upper probabilities, both conditional and unconditional.
In both cases, we write $\emph{coherent}$ meaning \emph{Williams-coherent} (shortly \emph{W-coherent}) \cite{williams_notes_2007}.
That is, we refer to the coherence concept developed by Williams, in the structure-free version studied in \cite{pelessoni_williams_2009}, where the set of events $\dset$ on which a lower/upper probability is defined may be arbitrary.
When $\dset$ is either finite or made of unconditional events only, W-coherence coincides with Walley's coherence \cite[Section 7.1.4 (b)]{walley_statistical_1991}, while being more general otherwise.
For the purposes of this paper, it is useful to think of W-coherence indirectly, not by means of its definition but through its characterisation given by the following \cite{pelessoni_williams_2009,williams_notes_2007}
\begin{theorem}[Envelope Theorem]
\label{thm:envelope_theorem}
Given a set of conditional events $\dset$, a lower probability $\lpr:\dset\rightarrow\rset$ is W-coherent on $\dset$ iff
\begin{eqnarray}
	\label{eq:envelope_def}
	\lpr(A|B)=\inf_{P\in\mset^*} P(A|B), \forall A|B\in\dset,
\end{eqnarray}
where $\mset^*$ is a set of (precise) probabilities.
\end{theorem}
We point out that the (precise) probabilities in $\mset^*$ are subjective or coherent in the sense of de Finetti.
This implies that $P(A|B)$ is defined also when ($B\in\dset$ and) $P(B)=0$.

There may be more sets $\mset^*$ satisfying Equation \eqref{eq:envelope_def}.
The largest such set (the set of all probabilities $P\geq\lpr$) is called the \emph{credal set} of $\lpr$ and denoted by $\mset$.
The infimum in \eqref{eq:envelope_def} is attained for $\mset^*=\mset$.
When $\dset$ is a set of unconditional events, $\mset$ is a convex set, while being not so in general.
The Envelope Theorem supports the interpretation of a coherent lower probability as a cautionary evaluation of a set of probabilities, one of which might be the `true' one, but there is uncertainty on which one.
With Walley's coherence, its `if' part obtains too, the `only if' implication not always in an infinite environment.

Coherent upper probabilities are similarly characterised by means of the Envelope Theorem, replacing \eqref{eq:envelope_def} with 
\begin{eqnarray*}
	\label{eq:envelope_def_con}
	\upr(A|B)=\sup_{P\in\mset^*} P(A|B), \forall A|B\in\dset.
\end{eqnarray*}
However, when considering simultaneously lower and upper probabilities, they will be \emph{conjugate}, i.e.
\begin{equation}
\label{eq:conjugacy}
\lpr(A|B)=1-\upr(\nega{A}|B).
\end{equation}
Equation \eqref{eq:conjugacy} lets us refer to lower (alternatively upper) probabilities only.
It is also important to recall that conjugate $\lpr$ and $\upr$ share the \emph{same} credal set $\mset$.

In this paper, we start with \emph{unconditional} lower probabilities ($\lpr(\cdot)$) and their conjugates ($\upr(\cdot)$).

Coherence implies that \cite[Sec. 2.7.4]{walley_statistical_1991}
\begin{equation}
\label{eq:coherence_necessary}
\text{if } A \Rightarrow B, \text{then } \lpr(A)\le\lpr(B), \upr(A)\le\upr(B) \text{  (\emph{monotonicity})}.
\end{equation}
The domain $\dset$ of $\lpr(\cdot)$, $\upr(\cdot)$ will often be $\asetpa$, the set of events logically dependent on a given partition $\prt$ (the powerset of $\prt$, in set theoretic language).

The following Lemma, whose simple proof employs \eqref{eq:conjugacy} and \eqref{eq:coherence_necessary}, will be useful later on.
\begin{lemma}
	\label{lem:pos_min_1}
	If $\lpr$, $\upr$ are coherent and conjugate on $\asetpa$ and $A, B\in\asetpa$, then
	\begin{align}
		\label{eq:pos_implies_1}
		\lpr(A\land B)>0 \text{ implies } \upr(\nega{A}\land B)<1;\\
		\label{eq:pos_implies_2}
		\upr(\nega{A}\land B)>0 \text{ implies } \lpr(A\land B)<1.
	\end{align}
\end{lemma}
A lower probability $\lpr$, coherent on $\asetpa$, is \emph{$2$-monotone} if $\lpr(A\vee B)+\lpr(A\wedge B)\geq\lpr(A)+\lpr(B)$, $\forall A,B\in\asetpa$.
Its conjugate $\upr$ is \emph{$2$-alternating}, meaning that $\upr(A\vee B)+\upr(A\wedge B)\leq\upr(A)+\upr(B)$, $\forall A,B\in\asetpa$.

From unconditional coherent lower/upper probability assessments on $\asetpa$, later on we shall consider their coherent extensions to $\asetpa\cup\asetpa|B$, where $B\in\asetpa\setminus\{\emptyset\}$ is given and
\begin{equation}
	\label{eq:alg_ext}
	\asetpa|B=\{A|B:A\in\asetpa\}.
\end{equation}
In other words, given a conditioning event $B$, we shall be interested in evaluating all events $A|B$, with $A$ varying in $\asetpa$.

In general, given a W-coherent lower probability $\lpr$ on an arbitrary set $\dset$ of conditional and/or unconditional events, $\lpr$ has a (W-)coherent extension, not necessarily unique, on any set of conditional events $\dset^\prime\supset\dset$.

The \emph{natural extension} $\lpe$ of $\lpr$ on $\dset^\prime$ is the \emph{least-committal} coherent extension of $\lpr$ to $\dset^\prime$, meaning that if $\lpq$ is a coherent extension of $\lpr$, then $\lpe\leq\lpq$ on $\dset^\prime$.
Further, $\lpe=\lpr$ on $\dset$ iff $\lpr$ is coherent \cite{pelessoni_williams_2009,walley_statistical_1991}.
$\lpe$ always exists with W-coherence, and is the preferred extension (when not the only one), since it incorporates only the information available from the starting assessment $\lpr$ on $\dset$ and nothing else.

Returning to the case of interest here, that $\lpr$ is coherent on $\asetpa$ and we look for its natural extension on $\asetpa\cup\asetpa|B$, we have by the Envelope Theorem that
\begin{equation}
\label{eq:nat_ext_expr}
\lpe(A|B)=\inf_{P\in\mset}P(A|B).
\end{equation}

Formula \eqref{eq:nat_ext_expr} specialises when $\lpr$ is $2$-monotone (and its conjugate $\upr$ is $2$-alternating), while $\lpr(B)>0$:
\begin{proposition}\emph{{(\cite[Theorem 7.2]{walley_coherent_1981}, \cite[Section 6.4.6]{walley_statistical_1991}).}}
\label{pro:2_monotone_properties}
If $\lpr$ is a coherent $2$-monotone lower probability on $\asetpa$ and $\upr$ is its conjugate, given $B\in\asetpa$ such that $\lpr(B)>0$, then, $\forall A\in\asetpa$,
\begin{align}
\label{eq:2_monotone_natural_extension_inf}
\lpe(A|B)&=\frac{\lpr(A\wedge B)}{\lpr(A\wedge B)+\upr(\nega{A}\wedge B)}\\
\label{eq:2_monotone_natural_extension_sup}
\upe(A|B)&=\frac{\upr(A\wedge B)}{\upr(A\wedge B)+\lpr(\nega{A}\wedge B)}.
\end{align}
$\lpe$ is $2$-monotone ($\upe$ is $2$-alternating) on $\asetpa|B$.
$\lpe$, $\upe$ are conjugate.
\end{proposition}

\subsection{Nearly-Linear Models}
\label{subsub:Nearly_Linear_Models}
Nearly-Linear models are a functionally simple instance of \emph{neighbourhood models},
given a (precise) probability $P_0$.
\begin{definition}
\label{def:Nearly_Linear_Models}
A \emph{Nearly-Linear Model} is a couple $(\lpr,\upr)$ of conjugate lower and upper probabilities on $\asetpa$,
where $\forall A\in\asetpa\setminus\{\impev,\Omega\}$
\begin{align}
\label{eq:Nearly_Linear_Models_1}
\lpr(A)&=\min\{\max\{bP_0 (A)+a,0\},1\},\\
\label{eq:Nearly_Linear_Models_2}
\upr(A)&=\max\{\min\{bP_0 (A)+c,1\},0\}
\end{align}
and $\lpr(\impev)=\upr(\impev)=0$, $\lpr(\Omega)=\upr(\Omega)=1$.

In Equations \eqref{eq:Nearly_Linear_Models_1}, \eqref{eq:Nearly_Linear_Models_2}, $P_0$ is an assigned probability on
$\asetpa$, while
\begin{equation}
\label{eq:Nearly_Linear_cofficients}
b>0,\ a\in\rset,\ c=1-(a+b).
\end{equation}
\end{definition}
NL models have been defined in \cite{corsato_nearly-linear_2019}, where their basic properties have been investigated.
It has been shown in \cite[Sec. 3.1]{corsato_nearly-linear_2019} that NL models are partitioned into three subfamilies, with varying consistency properties.
The coherent NL models are all the models in the VBM subfamily and some of the HBM (to be recalled next), while, within the third subfamily, $\lpr$ and $\upr$ are coherent iff the cardinality of $\prt$ is $2$ (therefore we neglect these latter models).
\begin{definition}
\label{def:Vertical_Barrier}
A \emph{Vertical Barrier Model} (VBM) is a NL model where \eqref{eq:Nearly_Linear_Models_1}, \eqref{eq:Nearly_Linear_Models_2}, \eqref{eq:Nearly_Linear_cofficients} specialise into
\begin{align}
\label{eq:Vertical_Barrier_Model_1}
\lpr(A)&=\max\{bP_0 (A)+a,0\},\ \forall A\in\asetpa\setminus\{\Omega\}, \ \lpr(\Omega)=1\\
\label{eq:Vertical_Barrier_Model_2}
\upr(A)&=\min\{bP_0 (A)+c,1\},\ \;\forall A\in\asetpa\setminus\{\impev\}, \ \;\upr(\impev)=0\\
\label{eq:Vertical_Barrier_Model_3}
0&\leq a+b\leq 1, a\leq 0
\end{align}
and $c$ is given by \eqref{eq:Nearly_Linear_cofficients} (hence $c\geq 0$).

In a \emph{Horizontal Barrier Model} (HBM) $\lpr$, $\upr$ are given by \eqref{eq:Nearly_Linear_Models_1}, \eqref{eq:Nearly_Linear_Models_2}, \eqref{eq:Nearly_Linear_cofficients} $\forall A\in\asetpa\setminus\{\impev,\Omega\}$, where $a$, $b$ satisfy the constraints
\begin{equation*}
\label{eq:HBM_coefficient_constraints}
a+b>1,\ b+2a\leq 1
\end{equation*}
(implying $a<0$, $b>1$, $c<0$).
\end{definition}
\begin{figure}[htbp!]
	\centering
	\includegraphics*[scale=0.25, trim = 1.8cm 0cm 0cm 0cm]{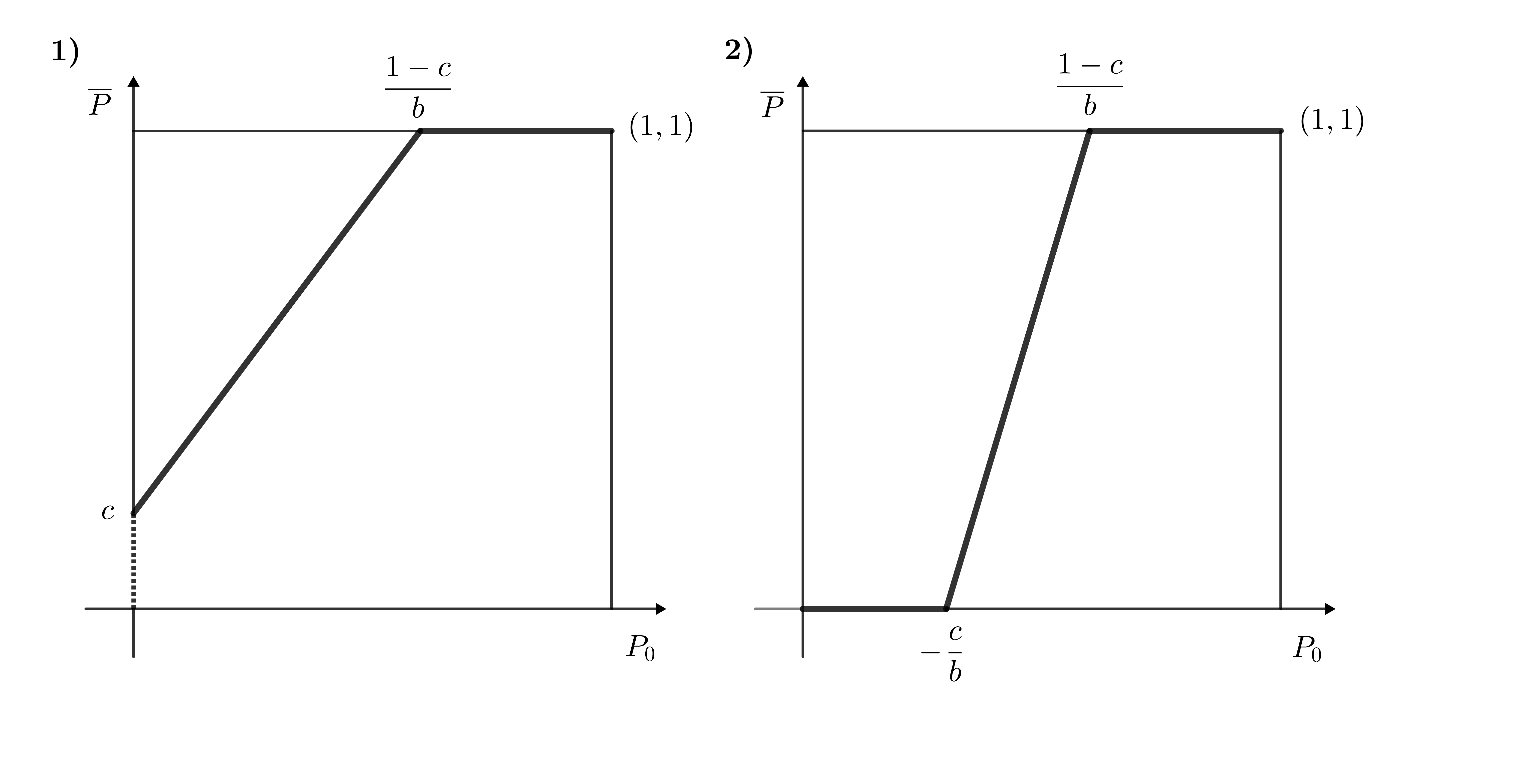}
	\caption{Plots of $\overline P$ (continuous bold line) against $P_0$: $1)$ in the VBM $2)$ in the HBM.}
	\label{fig:HBM_VBM_behaviour}
\end{figure}
In Figure \ref{fig:HBM_VBM_behaviour}, 1) we see the graph of a VBM $\upr$ against $P_0$ in the $(P_0,\upr)$ plane.
It originates a \emph{vertical barrier} (dotted) in the $\upr$ axis, meaning that no non-trivial event is given upper probability smaller than $c\ (<1)$.
By contrast, a HBM $\upr$ fixes a \emph{horizontal barrier} in the $P_0$-axis:
events whose $P_0$-probability is not larger than $-\frac{c}{b}$ are given upper probability $0$ in the HBM, see Figure \ref{fig:HBM_VBM_behaviour}, 2).

As for the coherence of VBMs and HBMs, we have that
\begin{proposition}\emph{(\cite[Propositions 4.2, 5.5]{corsato_nearly-linear_2019})}
\label{pro:NL_properties}
$\lpr$, $\upr$ are coherent and $2$-monotone, respectively $2$-al\-ter\-nat\-ing in any VBM; in a HBM they are so iff $\upr$ is subadditive (i.e. $\upr(A)+\upr(B)\geq\upr(A\vee B),\ \forall A,B\in\asetpa$).
\end{proposition}
\begin{remark}
\label{rem:marginal_HBM}
Thus, VBMs and (partly) HBMs ensure very good consistency properties.
We emphasise that subadditivity for HBMs is a considerably more restrictive requirement than it would seem at first glance. In fact, it is equivalent to more detailed and relatively uncommon conditions, as discussed in \cite[Section 5.2]{corsato_nearly-linear_2019}. This makes coherent HBMs somehow marginal within coherent NL models.
\end{remark}

A VBM generalises a number of well-known models.
Among them:
\begin{enumerate}
	\item[$\bullet$]
	if $a+b=0$, the \emph{vacuous lower/upper probability model} \cite[Sec. 2.9.1]{walley_statistical_1991}:
	\begin{align*}
	\lpr_V(A)&=0, \forall A\neq\Omega,\ \lpr_V(\Omega)=1, \\
	\upr_V(A)&=1, \forall A\neq\impev,\ \upr_V(\impev)=0;
	\end{align*}
	\item[$\bullet$]
	if $a=0$, $0<b<1$,
	the \emph{$\varepsilon$-contamination model} or linear-vacuous mixture model \cite[Sec. 2.9.2]{walley_statistical_1991}, here $b=1-\varepsilon$:
	\begin{align}
	\label{eq:vareps_formula_1}
	\lpr_{\varepsilon}(A)&=(1-\varepsilon) P_0 (A), \ \forall A\neq\Omega,\ \lpr_{\varepsilon}(\Omega)=1,\\
	\label{eq:vareps_formula_2}
	\upr_{\varepsilon}(A)&=(1-\varepsilon) P_0 (A)+\varepsilon,\ \forall A\neq\impev,\ \lpr_{\varepsilon}(\impev)=0;
	\end{align}
	\item[$\bullet$]
	if $b=1+\delta>1,\ a=-\delta<0$ (hence $c=0$), the \emph{Pari-Mutuel Model}
	\cite{montes_pari-mutuel_2019,pelessoni_inference_2010}, \cite[Sec. 2.9.3]{walley_statistical_1991}:
	\begin{equation*}
	\label{PMM_expression}
	\begin{aligned}[t]
		\lpr_{PMM}(A)&=\max\{(1+\delta)P_0(A)-\delta,0\},\\
	\upr_{PMM}(A)&=\min\{(1+\delta)P_0(A),1\};
	\end{aligned}
	\end{equation*}
	\item[$\bullet$]
	if $b=1$, $-1<a<0$ (hence $c=-a$), the \emph{Total Variation Model}
	\cite[Sec. 3]{herron_divisive_1997}, \cite[Sec. 3.2]{pelessoni_inference_2010}:\footnote{Note that $\lpr_{TVM}(A)\leq\upr_{TVM}(A), \forall A$, since $a<0$.}
	\begin{equation*}
	\label{TVM_expression}
	\begin{aligned}[t]
	\lpr_{TVM}(A)&=\max\{P_0 (A)+a,0\}\ \forall A\neq\Omega,\lpr_{TVM}(\Omega)=1,\\
	\upr_{TVM}(A)&=\min\{P_0 (A)-a,1\}\ \forall A\neq\impev,\ \upr_{TVM}(\impev)=0.
	\end{aligned}
	\end{equation*}
\end{enumerate}
VMBs may express a certain variety of consistent (and HBMs also inconsistent, although realistic) beliefs.
For these and other properties of these models, we refer to \cite{corsato_nearly-linear_2019, pelessoni_inference_2020}. 
\section{Conditioning Coherent Nearly-Linear Models}
\label{sec:NLM_conditioning}

This section is concerned with conditioning NL models via natural extension.
An alternative conditioning using the less familiar concept of regular extension is discussed in Section~\ref{subsec:cond_RE} in the Appendix.
%\subsection{Conditioning with the Natural Extension}
%\label{subsec:cond_NE}

Given a coherent NL model $(\lpr,\upr)$ on $\asetpa$
and an event $B\in\asetpa\setminus\{\impev\}$, we look for the natural extensions $\lpe(A|B)$, $\upe(A|B)$ of $\lpr$, $\upr$ respectively, for any $A\in\asetpa$.
In other words, $\lpr$, $\upr$ are extended on $\asetpa|B$ (defined in \eqref{eq:alg_ext}).

When $\lpr(B)=0$, we determine $\lpe$, $\upe$ quickly thanks to the next result, which applies more generally to coherent lower and upper probabilities and was proven in \cite[Proposition 3]{pelessoni_dilation_2020}.

\begin{proposition}
	\label{pro:natural_extension_case_zero_event}
	Let $\lpr:\dset\rightarrow\rset$ be a coherent lower probability on $\dset$, non-empty set of unconditional events, and $B\in\dset$, $B\neq\impev$ such that $\lpr(B)=0$. Then the natural extension $\lpe$ of $\lpr$ on $\dset\cup\{A_i|B\}_{i\in I}$, where $A_i\in\dset,\ \forall i\in I$, is given (by $\lpe(F)=\lpr(F)$, $\forall F\in\dset$ and) by
	\begin{equation}
		\label{eq:lower_natural_extension_zero_case_1}
		\lpe(A_i|B)=1 \text{ if } B\Rightarrow A_i,\ \lpe(A_i|B)=0 \mbox{ otherwise.}
	\end{equation}
	Correspondingly, the natural extension $\upe$ of the conjugate $\upr$ is given by $\upe(A_i|B)=1$ if $B\not\Rightarrow \nega{A_i}$, $\upe(A_i|B)=0$ if $B\Rightarrow \nega{A_i}$.
\end{proposition}
Note that Proposition \ref{pro:natural_extension_case_zero_event} points out an instance of \emph{vacuous} natural extension,
such that $\lpe$ ($\upe$) takes value $0$ (value $1$), and that in our framework $\dset=\asetpa$.

When $\lpr(B)>0$, $\lpe$, $\upe$ are instead determined from
\begin{proposition}\emph{(\cite[Proposition 4]{pelessoni_dilation_2020})}
\label{pro:NL_natural_extension}
Let $(\lpr,\upr)$ be a coherent NL model on $\asetpa$. For a given $B\in\asetpa$ such that $\lpr(B)>0$, we have that
\begin{align*}
%\label{eq:NL_natural_extension_2}
\lpe(A|B)=\begin{dcases}
1 \text{ iff } \upr(\nega{A}\wedge B)=0\\
\frac{bP_0(A\wedge B)+a}{bP_0(B)+1-b}\ (\in ]0,1[) \text{ iff }
\lpr(A\wedge B),\upr(\nega{A}\wedge B)>0\\
0 \text{ iff } \lpr(A\wedge B)=0
\end{dcases}\\
%\label{eq:NL_natural_extension_1}
\upe(A|B)=\begin{dcases}
	0 \text{ iff } \upr(A\wedge B)=0\\
	\frac{bP_0(A\wedge B)+c}{bP_0(B)+1-b}\ (\in ]0,1[) \text{ iff }
	\lpr(\nega{A}\wedge B),\upr(A\wedge B)>0\\
	1 \text{ iff } \lpr(\nega{A}\wedge B)=0
\end{dcases}
\end{align*}
\end{proposition}
\begin{remark}
\label{rem:diff_pos}
By Lemma \ref{lem:pos_min_1}, the hypothesis $\lpr(A\land B)>0$, $\upr(\nega{A}\land B)>0$ in the second entry for $\lpe(A|B)$ is equivalent to $\lpr(A\land B),\upr(\nega{A}\land B)\in ]0,1[$. This latter form was used in \cite[Proposition 4]{pelessoni_dilation_2020}.
The same remark applies to the second entry for $\upe(A|B)$ (exchange $A$ and $\nega{A}$ in Lemma \ref{lem:pos_min_1}).
\end{remark}

In the case of a VBM, the formulae for $\lpe$, $\upe$ in Proposition \ref{pro:NL_natural_extension} further specialise, as follows:
\begin{proposition}\emph{(\cite[Proposition 5]{pelessoni_dilation_2020})}
\label{pro:NE_VBM}
Let $(\lpr,\upr)$ be a VBM on $\asetpa$. For a given $B\in\asetpa$, with $\lpr(B)>0$, we have that
\begin{align}
\label{eq:NE_VBM_inf}
\lpe(A|B)&=\max\{b_B P_0(A|B)+a_B,0\}, \forall A\in\asetpa\setminus\{\Omega\}, \ \lpe(\Omega|B)=1 \\
\label{eq:NE_VBM_sup}
\upe(A|B)&=\min\{b_B P_0(A|B)+c_B,1\},\; \forall A\in\asetpa\setminus\{\impev\}, \ \; \upe(\impev|B)=0\\
\label{NE_VBM_coefficients}
a_B&=\frac{a}{bP_0(B)+1-b},\ b_B=\frac{bP_0(B)}{bP_0(B)+1-b},\ c_B=1-(a_B+b_B).
\end{align}
Moreover, it holds that
$b_B>0,\ a_B \leq 0,\ 0<a_B+b_B\leq 1$.
\end{proposition}
Proposition \ref{pro:NE_VBM} points out an important feature of a VBM: when all events in $\asetpa$ are conditioned on the same $B$, the resulting model is still a VBM.
(Note that this holds also when $\lpr(B)=0$:
here Proposition \ref{pro:NE_VBM} does not apply, but from Proposition~\ref{pro:natural_extension_case_zero_event} we obtain the vacuous lower/upper probabilities, a special VBM.)

We express this property saying that a VBM is \emph{stable under conditioning}.

\begin{remark}
\label{rem:how_to_extend}
As an interesting follow-up of stability of VBMs and of results in \cite{pelessoni_inference_2020}, we may obtain the natural extension $\lpe$ to a conditional gamble $X|B$ when $\lpr(B)>0$, starting from a VBM on $\asetpa$.
Recall that a gamble is a bounded random number, and that the natural extension of $\lpr$ to a gamble is again the least-committal coherent extension of $\lpr$, cf. \cite{pelessoni_williams_2009,troffaes_lower_2014,walley_statistical_1991}.
We require that event $(X=x_i)$ belongs to $\asetpa$, for any possible value $x_i$ of $X$.
\begin{figure}[htbp!]
	\centering
	\includegraphics*[scale=0.550, trim = 8.0cm 15.5cm 13cm 5.5cm]{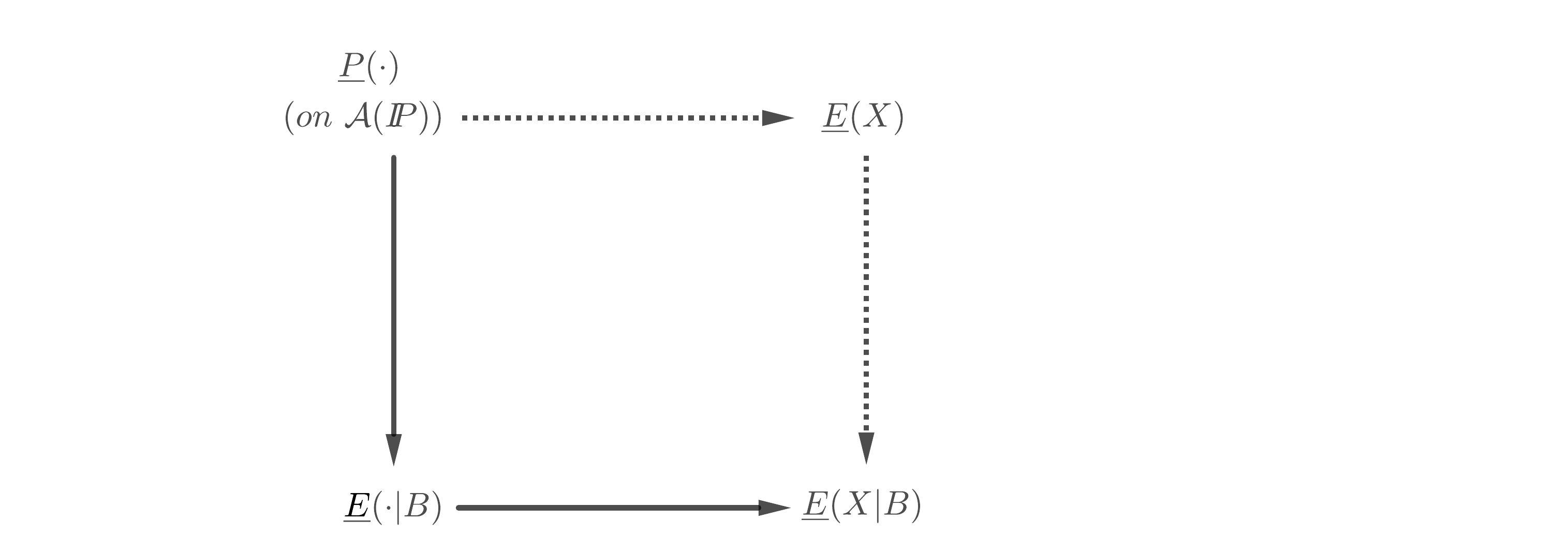}
	\caption{Natural extension of a VBM on $X|B$.}
	\label{fig:VBM_natural_extension}
\end{figure}
To compute $\lpe(X|B)$, we may first determine $\lpe$ on $\asetpa|B$ by Proposition \ref{pro:NE_VBM}.
To obtain $\lpe(X|B)$ we apply then \cite[Proposition 4.3]{pelessoni_inference_2020}, where $X$ is replaced by $X|B$ everywhere (this proposition is stated in \cite{pelessoni_inference_2020} in unconditional terms, i.e. conditional on $\Omega$, but it applies whenever it involves a single conditioning event $B$).
This two-step procedure corresponds to following the continuous route in Figure \ref{fig:VBM_natural_extension}, which appears simpler than the alternative dotted route.
\end{remark}

We mention that the main VBM submodels are stable too \cite{pelessoni_dilation_2020}.
This means that the natural extension on $\asetpa |B$ of a Pari-Mutuel Model or a Total Variation Model is again, respectively, a Pari-Mutuel Model or a Total Variation Model.
The same property holds with $\varepsilon$-contamination models and (trivially) linear-vacuous models.\footnote{
	The stability properties of these submodels and their natural extensions have been derived independently also in
	\cite[Sections 4.2, 5.2]{montes_neighbourhood_I}, \cite[Section 2.2]{montes_neighbourhood_II}.
} 
As shown in \cite[Example 1]{pelessoni_dilation_2020}, HBMs are instead generally \emph{not stable} under conditioning.
We refer to \cite[Section 3]{pelessoni_dilation_2020} for a detailed discussion of stability.

\section{Dilation and Nearly Linear Models}
\label{sec:dilation_NL}
\subsection{What is dilation?}
\label{subsec:what_dilation}
Let us recall the definition of weak and strict dilation in general \cite{herron_divisive_1997,pelessoni_dilation_2020}, when coherent $\lpr(A)$, $\upr(A)$ are given and conditioning $A$ on each event of a partition is performed by means of the natural extensions $\lpe$, $\upe$.
\begin{definition}
	\label{def:weakly_dilation}
	Given an event $A$ and a partition $\bset$, say that $\bset$ \emph{weakly dilates} $A$ iff
	\begin{equation}
		\label{eq:weakly_dilation}
		\lpe(A|B)\le\lpr(A)\le\upr(A)\le\upe(A|B), \forall B\in\bset.
	\end{equation}
	
	Weak dilation is \emph{trivial} if \eqref{eq:weakly_dilation} holds with its three inequalities being equalities, $\forall B\in\bset$, is \emph{non-trivial} otherwise.
	
	$\bset$ \emph{strictly dilates} $A$ if the outer inequalities in \eqref{eq:weakly_dilation} are strict for any $B\in\bset$.
\end{definition}
Note that dilation implies the weaker phenomenon of \emph{imprecision increase}, occuring when $\upe(A|B)-\lpe(A|B)\geq\upr(A)-\lpr(A)$, $\forall B\in\bset$.
Imprecision increase is an already puzzling phenomenon:
it ensures that the conditional evaluation is more imprecise than the unconditional one, or at least equally imprecise, \emph{whatever happens}, i.e. whatever conditioning event $B$ in partition $\bset$ turns out to be true.

A first question is whether some assumptions can be introduced, to rule out instances where it is already known whether dilation occurs or not. Recall for this that
\begin{definition}
	\label{def:extreme_event}
	An event $A$ is \emph{logically independent} of a partition $\bset$ if $B\not\Rightarrow A$ and $B\not\Rightarrow\nega{A}$, $\forall B\in\bset$.
	
	An event $A$ is called \emph{extreme} if $\lpr(A)=\upr(A)=0$ or $\lpr(A)=\upr(A)=1$.
\end{definition}
Then, the results discussed in detail in Section \ref{sec:dilation_generalities} in the Appendix let us conclude that

\emph{the investigation of dilation can be restricted to events that are logically independent of $\bset$ and, whenever conditioning is performed by means of the natural extension, non-extreme}.

\subsection{Dilation with Nearly Linear Models - Basic Assumptions}
\label{subsec:dilation_with_NL}
Let us consider now a NL model.
This adds some more specific hypotheses:
$\lpr$ and its conjugate $\upr$ are defined on $\asetpa$, and as for partition $\bset$, we assume that
\begin{equation}
\label{eq:partitions}
\bset\subset\asetpa\setminus\{\emptyset\}.
\end{equation} 
%\subsection{Basic Assumptions}
%\label{sec:basic_assumptions}
The conclusions of Section \ref{subsec:what_dilation} let us also assume from now onwards that the following $(A1)$, $(A2)$ apply. We also require $(A3)$ to hold, which rules out a degenerate situation.\footnote{
Precisely, suppose $(A3)$ does not apply, so that $\lpr(A)=0, \upr(A)=1$.
Then, Equation~\eqref{eq:weakly_dilation} ensures that either $(a)$ dilation does not occur, or $(b)$ it is degenerate, in the sense that also $\lpe(A|B)=0, \upe(A|B)=1, \forall B\in\bset$.
While this is true in general, with coherent NL models and assuming $(A1)$, only the degenerate case $(b)$ occurs.
This follows from Proposition~\ref {pro:natural_extension_case_zero_event} when $\lpr(B)=0$.
Instead, when $\lpr(B)>0$, it ensues from Proposition~\ref{pro:NL_natural_extension},
noting that $\lpr(A \wedge B)=0$, hence $\lpe(A|B)=0$, and that $\lpr(\nega{A}\wedge B)=1-\upr(A\vee\nega{B})=0$ because $1\geq\upr(A\vee\nega{B})\geq\upr(A)=1$,
hence $\upe(A|B)=1$.
}
\begin{enumerate}
\item[\textbf{(A1)}] $A$ is logically independent of $\bset$.
\item[\textbf{(A2)}] $A$ is non-extreme.
\item[\textbf{(A3)}] $(\lpr(A),\upr(A))\neq(0,1)$.
\end{enumerate}

\subsection{The Role of Logical Independence}
\label{sec:role_log_ind}
Preliminarily, it is interesting to point out that, essentially because of the assumption \eqref{eq:partitions}, we may state that the logical independence condition $(A1)$ does not apply for a number of choices of $A$ and $\bset$ (where, consequently, dilation does or does not occur, according to whether, respectively, Corollary \ref{cor:1-semidependence} or Proposition \ref{pro:indep_dilation} applies).

To see this, we write $|A|$ for the \emph{cardinality} of $A$, meaning the cardinality of the set $\{\omega\in\prt: \omega\Rightarrow A\}$.
Similarly, $|\bset|$ is the cardinality of $\{B: B\in\bset\}$.
Then, the following lemma holds.
\begin{lemma}
	\label{lem:cardinality}
	Given $A\in\asetpa$, $\bset\subset\asetpa\setminus\{\emptyset\}$, condition 
	\begin{equation*}
	\label{eq:cardinality}
	2\leq |\bset|\leq\min\{|A|, |\nega{A}|\}
	\end{equation*}
	is necessary for $A$ to be logically independent of $\bset$.
\end{lemma}
\begin{proof}
	Necessarily $|\bset|\geq 2$, or else $\bset$ is the trivial partition $\bset=\{\Omega\}$.
	
	Further, for any $B\in\bset$ there must exist two distinct $\omega, \omega^\prime\in\prt$ such that $\omega\Rightarrow A\land B$, $\omega^\prime\Rightarrow\nega{A}\land B$.
	In fact, this guarantees that $A\land B\neq\emptyset$, $\nega{A}\land B\neq\emptyset$, equivalent to the logical independence conditions $B\not\Rightarrow\nega{A}$, $B\not\Rightarrow A$, respectively.
	Therefore, $|\bset|\leq\min\{|A|, |\nega{A}|\}$.
	The bound is tight, being reached when for any $B\in\bset$ there is exactly one atom of $\prt$ implying both $B$ and the event between $A$ and $\nega{A}$ with smaller cardinality (if the cardinality is the same, two such atoms are needed for each $B$).
\end{proof}
The following facts are easy follow-ups of Lemma \ref{lem:cardinality}:
\begin{enumerate}
\item[(a)] When $\prt$ is infinite, $A$ is not logically independent of $\bset$ if $\bset$ is infinite and $A$ is either finite or cofinite.
\item[(b)] When $\prt$ is finite and made up of $n$ atoms, $|A|$ constrains the type and number of the partitions $A$ may be logically independent of. For instance:
\begin{enumerate}
	\item[-] if $|A|=1$ or $|A|=n-1$ (i.e. either $A$ or $\nega{A}$ is an atom of $\prt$), $A$ is logically independent of no partition $\bset\subset\asetpa$;
	\item[-] if $|A|=2$ or $|A|=n-2$ (and $n\geq 4$), $A$ is logically independent of $\bset$ only if $\bset$ is binary. There are $4(2^{n-3}-1)$ such partitions that are admissible;
	\item[-] if $|A|=3$ or $|A|=n-3$ (and $n\geq 6$), $\bset$ can be either binary (there are $3(2^{n-3})$ such admissible partitions) or ternary (with $3!\binom{n-3}{3}\cdot 3!\binom{n-4}{2}=3(n-3)(n-4)^2 (n-5)^2$ distinct admissible partitions).
\end{enumerate}
\end{enumerate}
We see from $(a)$ and $(b)$ above that dilation can never affect \emph{all} non-trivial events of $\asetpa$.

As a general qualitative rule, the higher the cardinality of $\bset$ given $\prt$, the more logical independence, hence dilation, is restricted to fewer or in the limit no events $A$.
When $\prt$ is finite and formed by $n$ atoms, it suffices to prevent logical independence that $|\bset|>\lfloor{\frac{n}{2}}\rfloor$, since then necessarily at least one $B\in\bset$ is an atom of $\prt$.
This makes any $A$ semidependent or dependent on $\bset$.

\subsection{Results on Dilation}
\label{sec:results_dilation}
Dilation for coherent NL models is characterised in Proposition \ref{pro:dilation_NLM}, relying upon the next preliminary lemma.
\begin{lemma}
	\label{lem:partial_dilation}
	Given a coherent NL model $(\lpr,\upr)$, let $A, B\in\asetpa$,
	$\lpr(B)>0$.
	\begin{enumerate}
		\item[(i)]
		If $\lpr(A\wedge B), \upr(\nega{A}\wedge B)>0$, then
		\begin{equation*}
			%\lpe(A|B)\leq\lpr(A) \text{ iff either } P_0(\nega{B})=0 \text{ or } \lpr(A)\le P_0(A|\nega{B}).
			\lpe(A|B)\leq\lpr(A) \text{ iff } P_0(A \land \nega{B})\geq P_0(\nega{B})\lpr(A).
		\end{equation*}
		\item[(ii)]
		If $\lpr(\nega{A}\wedge B), \upr(A\wedge B)>0$, then
		\begin{equation*}
			%\upe(A|B)\geq\upr(A) \text{ iff either } P_0(\nega{B})=0 \text{ or } \upr(A)\ge P_0(A|\nega{B}).
			\upe(A|B)\geq\upr(A) \text{ iff } P_0(A \land \nega{B})\leq P_0(\nega{B})\upr(A).
		\end{equation*}
		\item[(iii)]
			If $\lpr(A\land B)=0$, then $\lpe(A|B)\leq\lpr(A)$.
		\item[(iv)]
			If $\lpr(\nega{A}\land B)=0$, then $\upe(A|B)\geq\upr(A)$.
	\end{enumerate}
\end{lemma}
\begin{proof}
	$(i)$ and $(ii)$ have been proven in \cite[Proposition 6]{pelessoni_dilation_2020};
	$(iii)$ follows from \eqref{eq:2_monotone_natural_extension_inf}, since $\lpe(A|B)=0$,
	$(iv)$ from \eqref{eq:2_monotone_natural_extension_sup}, since $\upe(A|B)=1$.
\end{proof}

\begin{proposition}
	\label{pro:dilation_NLM}
	Given a coherent NL model $(\lpr,\upr)$ on $\asetpa$, a partition $\bset\subset\asetpa\setminus\{\impev\}$ and $A\in\asetpa$ non-extreme and logically independent of $\bset$,
	\begin{enumerate}
	\item[1)]
	Equation \eqref{eq:weakly_dilation} holds for $B\in\bset$, $\lpr(B)>0$ iff $(a)$ and $(b)$ both apply:
	\begin{enumerate}
		\item[(a)]
		one of the following two conditions holds:
		\begin{enumerate}
			\item[(a1)]
			$\lpr(A\wedge B), \upr(\nega{A}\wedge B)>0$ and $P_0(A \land \nega{B})\geq P_0(\nega{B})\lpr(A)$, or 
			\item[(a2)]
			$\lpr(A\land B)=0$;
		\end{enumerate}
		\item[(b)]
		one of the following two conditions holds:
		\begin{enumerate}
			\item[(b1)]
			$\lpr(\nega{A}\wedge B), \upr(A\wedge B)>0$ and $P_0(A \land \nega{B})\leq P_0(\nega{B})\upr(A)$, or
			\item[(b2)]
			$\lpr(\nega{A}\land B)=0$.
		\end{enumerate}
	\end{enumerate}
	\item[2)]
		Weak dilation occurs (w.r.t. $A$, $\bset$) iff $(a)$ and $(b)$ both apply $\forall B\in\bset$ such that $\lpr(B)>0$.
	\end{enumerate}
\end{proposition}
\begin{proof}
\begin{enumerate}
\item[1)]
If $(a)$ and $(b)$ both hold, \eqref{eq:weakly_dilation} follows straightforwardly from Lemma \ref{lem:partial_dilation}.
	
We prove that $(a)$ and $(b)$ are necessary for \eqref{eq:weakly_dilation}. For this, let \eqref{eq:weakly_dilation} hold. Consider Proposition \ref{pro:NL_natural_extension} and note that
no couples of possible values of $\lpe(A|B)$, $\upe(A|B)$, other than those implied by $(a)$ and $(b)$, are compatible with \eqref{eq:weakly_dilation} and the hypotheses.
In fact, if $\lpr(\nega{A}\land B)=\upr(\nega{A}\land B)=0$ or $\upr(A\land B)=\lpr(A\land B)=0$, \eqref{eq:weakly_dilation} implies $\lpr(A)=\upr(A)=1$ or $\lpr(A)=\upr(A)=0$ respectively. Hence $A$ would be extreme, against the hypotheses. In the remaining discarded cases, Proposition~\ref{pro:NL_natural_extension}  implies $\upe(A|B)<\lpe(A|B)$, a contradiction.
\item[2)]
We need not consider those $B\in\bset$ such that $P(B)=0$, if any. In fact, they do not hinder dilation: logical independence of $A$ from $\bset$ implies for them $\lpe(A|B)=0$, $\upe(A|B)=1$ by Proposition \ref{pro:natural_extension_case_zero_event}. For the others apply $1)$.
\end{enumerate}
\end{proof}
We apply Proposition \ref{pro:dilation_NLM} in an example.
\begin{example}
	\label{ex:dilation}
A VBM is obtained using \eqref{eq:Nearly_Linear_cofficients}, \eqref{eq:Vertical_Barrier_Model_1}, \eqref{eq:Vertical_Barrier_Model_2} from a given probability $P_0$ on a six-atom partition $\prt=\{\omega_1,\ldots,\omega_6\}$. Take $A=\omega_2\lor\omega_4\lor\omega_5$ and partition $\bset=\{B_1, B_2, B_3\}$, with $B_i=\omega_{2i-1}\lor\omega_{2i}$, $i=1,2,3$.
% $B_1=\omega_1\lor\omega_2$, $B_2=\omega_3\lor\omega_4$, $B_3=\omega_5\lor\omega_6$.
\begin{table}[htbp!]
	\setlength{\tabcolsep}{10pt}
	\begin{center}
		\caption{Data for Example \ref{ex:dilation}}
		\begin{tabular}{c||c|c||c|c||c|c||c|c|c}
			\label{tab:example_coarsening}
			&\multicolumn{2}{|c||}{\,}&\multicolumn{2}{|c||}{\,}&\multicolumn{2}{|c||}{\,}&\multicolumn{2}{|c}{\,}&\\
			[-1em] & \multicolumn{2}{|c||}{$B_1$} & \multicolumn{2}{|c||}{$B_2$}   & \multicolumn{2}{|c||}{$B_3$}    &  \multicolumn{3}{|c}{\,}\\
			\hline &&&&&&&&\\[-1em]
			&&&&&&&&\\[-1em] $\prt$ & $\omega_1$ & $\omega_2$ & $\omega_3$ & $\omega_4$ & $\omega_5$ & $\omega_6$ & $b$ & $a$ & $c$\\
			&&&&&&&&\\[-1em]\hline &&&&&&&&\\[-1em]
			$P_0$ & 0.1 & 0.2 & 0.1 & 0.1 & 0.25 & 0.25 & 1.1 & -0.2 & 0.1
		\end{tabular}
		\label{tab:HBM_conditioning}
	\end{center}
\end{table}
Using the data in Table~\ref{tab:example_coarsening}, we have
$\lpr(A)=0.405,\, \upr(A)=0.705$.

Applying Proposition \ref{pro:dilation_NLM} to $A$ and $\bset$, we obtain:
\begin{itemize}
	\item[$\bullet$]
	With $B_1$, using $(a1)$ and $(b2)$:
	\begin{itemize}
	\item[]
	$\lpr(A\land B_1)=\lpr(\omega_2)=0.02$, $\upr(\nega{A}\land B_1)=\upr(\omega_1)=0.21$, and
	\item[]
	$P_0(A\land\nega{B_1})=P_0(\omega_4\lor\omega_5)=0.35>P_0(\nega{B_1})\lpr(A)=0.2835$, while
	\item[]
	$\lpr(\nega{A}\land B_1)=\lpr(\omega_1)=0$.
	\end{itemize}
	\item[$\bullet$]
	With $B_2$, $(a2)$ and $(b2)$ obtain:
	\begin{itemize}
	\item[]
	$\lpr(A\land B_2)=\lpr(\nega{A}\land B_2)=0$.
	\end{itemize}
	\item[$\bullet$]
	With $B_3$, cases $(a1)$ and $(b1)$ apply.
	In fact, it can be checked that
	\begin{itemize}
	\item[]
	$\lpr(A\land B_3)=0.075$, $\upr(\nega{A}\land B_3)=0.375$ and
	\item[]
	$P_0(A\land\nega{B_3})=0.3>P_0(\nega{B_3})\lpr(A)=0.2025$, while
	\item[]
	$\lpr(\nega{A}\land B_3)=0.075$, $\upr(A\land B_3)=0.375$, and
	\item[]
	$P_0(A\land\nega{B_3})=0.3<P_0(\nega{B_3})\upr(A)=0.3525$.
	\end{itemize}
\end{itemize}
Thus, $\bset$ dilates $A$.
\end{example}
\begin{remark}
\label{rem:dil_suff}
\emph{Comments on Proposition \ref{pro:dilation_NLM}.}

Concerning the \emph{interpretation} of the dilation conditions in Proposition~\ref{pro:dilation_NLM},
we observe that:
\begin{itemize}
\item[(a)]
The extreme assignments $\lpr(A\land B)=0$, $\lpr(\nega{A}\land B)=0$ favour dilation.
\item[(b)]
In their absence and assuming $P_0(\nega{B})>0$, $\forall B\in\bset$ (this rules out $P_0$-probabilities concentrated on a single atom of partition $\bset$), dilation occurs iff
\begin{equation}
	\label{eq:dilation_rules_out}
	\lpr(A)\leq P_0(A|\nega{B})\leq\upr(A), \forall B\in\bset \text{ s.t. } \lpr(B)>0.
\end{equation}
\item[(c)]
Let the assumptions in $(b)$ apply to a VBM. For VBMs, we know that $\lpr(A)\leq P_0(A)\leq\upr(A)$.
Thus, \eqref{eq:dilation_rules_out} is certainly verified if $A$ is \emph{$P_0$-non-correlated} with any $\nega{B}$,
since then $P_0(A|\nega{B})=P_0(A)$.
More generally, dilation is ensured by \eqref{eq:dilation_rules_out} if, for any $B$, $P_0(A|\nega{B})$ does not deviate too much from $P_0(A)$, since both must belong to the interval $[\lpr(A), \upr(A)]$.
\end{itemize} 
\end{remark}

Proposition \ref{pro:dilation_NLM} $2)$ offers a unified characterisation of dilation for all NL models. It extends previous results for special VBM submodels stated in simplified but heterogeneous forms \cite{herron_divisive_1997,pelessoni_inference_2010,seidenfeld_dilation_1993}.
To give an idea of this, we derive from Proposition \ref{pro:dilation_NLM} the characterisation for (weak) dilation with $\varepsilon$-contamination models in \cite[Proposition~1]{herron_divisive_1997}.
\begin{proposition}
\label{dilation-epsilon}
Given an $\varepsilon$-contamination model $(\lpr_\varepsilon,\upr_\varepsilon)$ (with $\varepsilon>0$), a partition $\bset\subset\asetpa\setminus\{\impev\}$ and $A\in\asetpa$ non-extreme and logically independent of $\bset$, suppose that
\begin{equation}
	\label{eq:epsilon_condition}
	P_0(A\land B)>0, P_0(\nega{A}\land B)>0, \forall B\in\bset.
\end{equation}
%Then, $\forall B\in\bset$ condition \eqref{eq:single_dilation} boils down to
Then, $\bset$ dilates (weakly) $A$ iff, $\forall B\in\bset$,
\begin{equation}
	\label{eq:epsilon_boil_down}
	-\varepsilon P_0(\nega{A})P_0(\nega{B})\leq P_0(A\land B)- P_0(A) P_0(B)\leq \varepsilon P_0(A)P_0(\nega{B}).
\end{equation}
\end{proposition}
\begin{proof}
Recalling \eqref{eq:vareps_formula_1}, \eqref{eq:vareps_formula_2} and since $A$ is non-extreme,
we easily deduce that conditions \eqref{eq:epsilon_condition} are equivalent, both when $A'=A$ and when $A'=\nega{A}$, to
$\lpr_\varepsilon (A' \land B)=(1-\varepsilon)P_0(A' \land B)>0$, hence $\upr_\varepsilon (A' \land B)=(1-\varepsilon)P_0(A' \land B)+\varepsilon>0$ as well, $\forall B\in\bset$.

Therefore, by Proposition \ref{pro:dilation_NLM} (apply $(a1)$ and $(b1)$ for any $B$), $\bset$ dilates (weakly) $A$ iff, $\forall B\in\bset$,
\begin{equation}
	\label{eq:epsilon_boil_down_2}
	(1-\varepsilon)P_0(A)P_0(\nega{B})\leq P_0(A \land \nega{B})\leq (1-\varepsilon)P_0(A)P_0(\nega{B})+\varepsilon P_0(\nega{B})
\end{equation}
To show that \eqref{eq:epsilon_boil_down} and \eqref{eq:epsilon_boil_down_2} are equivalent,
subtract $P_0(A)P_0(\nega{B})$ in \eqref{eq:epsilon_boil_down_2}:
\begin{equation*}
-\varepsilon P_0(A)P_0(\nega{B})\leq P_0(A \land \nega{B})-P_0(A)P_0(\nega{B})\leq \varepsilon P_0(\nega{A})P_0(\nega{B}).
\end{equation*}
Now substitute $P_0(A \land \nega{B})-P_0(A)P_0(\nega{B})=P_0(A) P_0(B)-P_0(A\land B)$ in the last expression, getting
\begin{equation*}
-\varepsilon P_0(A)P_0(\nega{B})\leq P_0(A) P_0(B)-P_0(A \land B)\leq \varepsilon P_0(\nega{A})P_0(\nega{B}),
\end{equation*}
equivalent to \eqref{eq:epsilon_boil_down}. 
\end{proof}
Strictly speaking, condition \eqref{eq:epsilon_boil_down} as derived in Proposition \ref{dilation-epsilon} concerns weak dilation, while \cite{herron_divisive_1997} focuses on strict dilation, where the inequalities in \eqref{eq:epsilon_boil_down} are strict. The argument for proving the equivalence with the strict version of \eqref{eq:epsilon_boil_down} is analogous.

Remark \ref{rem:dil_suff} obviously applies to the $\varepsilon$-contamination model as a special VBM.
Its peculiarity (recall \eqref{eq:vareps_formula_1} and \eqref{eq:vareps_formula_2}) is that, concerning $(a)$ in the remark, $\lpr_{\varepsilon}(A\wedge B)=\lpr_{\varepsilon}(\nega{A}\wedge B)=0$ implies $P_0(A)=0$;
as for $(c)$, $\upr_{\varepsilon}(A)-\lpr_{\varepsilon}(A)=\varepsilon$,
so that $P_0(A|\nega{B})$ cannot deviate from $P_0(A)$ by more than $\varepsilon$ for dilation to occur (while both belonging to $[\lpr_{\varepsilon}(A),\upr_{\varepsilon}(A)])$.

\section{Properties of Dilation for VBMs}
\label{sec:properties_dilation_VBM}
In this section we focus on VBMs to discuss various properties of dilation.
The reasons for laying aside HBMs are that unlike VBMs these models are very often not coherent (cf. Remark \ref{rem:marginal_HBM}), and even when they are, their conditional models are generally not of the same type, being not stable.

On the other hand, we shall see that our results extend previous ones, obtained for special cases of VBMs.
\subsection{The Coarsening Property}
\label{sec:coarsening}

When a partition dilates an event, it may be interesting to know whether this also happens with a different partition.
The coarsening property, introduced in \cite{herron_divisive_1997}, is related to this problem.
\begin{definition}
\label{def:coarsening}
Say that an imprecise probability model satisfies the \emph{coarsening property} if,
whenever a partition $\bset$ with three or more atoms dilates an event $A$, there is a (non-trivial) partition $\bset^\prime\neq\bset$, coarser than $\bset$, that also dilates $A$.
\end{definition}
%The model satisfies the \emph{binary dilation property} in the special case that partition $\bset^\prime$ is binary.
When the model is a VBM, we already know from \cite{herron_divisive_1997} that it has the coarsening property if it is an $\varepsilon$-contamination model, $\bset$ is finite, and the relevant events have no extreme probabilities.

Our goal in this section is to detect conditions for a VBM to satisfy the coarsening property, under assumptions generally weaker than those in \cite{herron_divisive_1997}.

The first step is the following
\begin{definition}
\label{def:elle}
Given a VBM on $\asetpa$, define for any $A, B\in\asetpa$
\begin{equation}
\label{eq:elle}
\lset(A,B)=P_0(A\land B)-b P_0(A) P_0(B)-a P_0(B).
\end{equation}
\end{definition}
\begin{proposition}
%[Properties of $\lset(A,B)$]
\label{elle_properties}
If $B=\vee_{i=1}^{n} C_i$, $C_i\land C_j=\emptyset$ if $i\neq j$, $C_i\in\asetpa\ (i=1,\ldots,n)$, then
\begin{equation}
	\label{eq:sum_L}
\lset(A,B)=\sum_{i=1}^{n}\lset(A,C_i).
\end{equation}
\end{proposition}
\begin{proof}
\begin{equation*}
\begin{aligned}[t]
	\sum_{i=1}^{n}\lset(A,C_i)&=&\sum_{i=1}^{n} P_0(A\land C_i)-b P_0(A)\sum_{i=1}^{n} P_0(C_i)-a \sum_{i=1}^{n} P_0(C_i)\\
	&=&P_0(A\land B)-b P_0(A) P_0(B)-a P_0(B)=\lset(A,B).
\end{aligned}
\end{equation*}
\end{proof}
\begin{corollary}
\label{cor:elle_properties}
Let $(\lpr,\upr)$ be a VBM on $\asetpa$.
% (defined as in \eqref{eq:Vertical_Barrier_Model_1}).
If $\bset=\{B_1,\ldots,B_n\}\subset\asetpa$ is a partition of the sure event $\Omega$, then
\begin{equation}
\label{eq:elle_inequality}
\sum_{i=1}^{n}\lset(A,B_i)\geq P_0(A)-\lpr(A)\geq 0.
\end{equation}
\end{corollary}
\begin{proof}
Put $C_i=B_i$ in Equation \eqref{eq:sum_L}. By \eqref{eq:elle} and since $b P_0(A)+a\leq\lpr(A)$ by \eqref{eq:Vertical_Barrier_Model_1}, while $P_0(A)\geq\lpr(A)$ in a VBM, we obtain
\begin{equation*}
\sum_{i=1}^{n}\lset(A,B_i)=\lset(A,\Omega)=P_0(A)-b P_0(A)-a\geq P_0(A)-\lpr(A)\geq 0.
\end{equation*}
\end{proof}
The significance of the numbers $\lset(\cdot,\cdot)$ is due to their capability of detecting in a number of instances whether $\lpe(A|B)\leq\lpr(A)$ and $\upe(A|B)\geq\upr(A)$, as we will now show.
\begin{proposition}
	\label{pro:elle_use}
	Given a VBM on $\asetpa$, let $\bset\subset\asetpa\setminus\{\emptyset\}$ be a partition.
	Take $A\in\asetpa$ and $B\in\bset$ such that $\lpr(B)>0$.
	\begin{enumerate}
		\item[(a)]
		If one of the following two conditions applies
		\begin{itemize}
			\item[(a1)]
			$\lpr(A)=0$,
			\item[(a2)]
			$\lpr(A\wedge B), \upr(\nega{A}\wedge B)>0$,
		\end{itemize}
	then
	\begin{equation}
		\label{eq:L_nonneg_1}
		\lpe(A|B)\leq\lpr(A) \text{ iff } \lset(A,\nega{B})\geq 0.
	\end{equation} 
	\item[(b)]
		If one of the following two conditions applies
		\begin{itemize}
			\item[(b1)]
			$\upr(A)=1$,
			\item[(b2)]
			$\lpr(\nega{A}\wedge B), \upr(A\wedge B)>0$,
		\end{itemize}
	then
	\begin{equation}
		\label{eq:L_nonneg_2}
		\upe(A|B)\geq\upr(A) \text{ iff } \lset(\nega{A},\nega{B})\geq 0.
	\end{equation} 
	\end{enumerate}
\end{proposition}
\begin{proof}
	\begin{enumerate}
		\item[(a)]
		If $(a1)$ holds, $\lpr(A\land B)=0$ and by \eqref{eq:2_monotone_natural_extension_inf} $\lpe(A|B)=0\leq\lpr(A)$.
		On the other hand, when $\lpr(A)=0$, by \eqref{eq:Vertical_Barrier_Model_1} we have that $b P_0(A)+a\leq 0$, and therefore $	\lset(A,\nega{B})=P_0(A\land\nega{B})+P_0(\nega{B})(-bP_0(A)-a)\geq 0$.
			
		Let now $(a2)$ hold. Since $\lpr(A\land B)>0$, also $\lpr(A)>0$, while $\lpr(A)<1$ by $(A2)$, Section~\ref{subsec:dilation_with_NL}.
		Hence, $\lpr(A)=bP_0(A)+a$ by \eqref{eq:Nearly_Linear_Models_1}.
		Using this when applying Lemma \ref{lem:partial_dilation} $(i)$ and Definition~\ref{def:elle}, we have that
		$\lpe(A|B)\leq\lpr(A)$ iff $P_0(A\land \nega{B})\geq P_0(\nega{B})(bP_0(A)+a)$ iff $\mathcal{L}(A,\nega{B})\geq 0$.
		\item[(b)] Follows from $(a)$, replacing $A$ with $\nega{A}$ and exploiting conjugacy.
	\end{enumerate}
\end{proof}
The next theorem concerns the coarsening property.
\begin{theorem}
	\label{thm:coarsening_basic}
	Given a VBM $(\lpr,\upr)$ on $\asetpa$, let $\bset\subset\asetpa\setminus\{\emptyset\}$ be a finite partition ($|\bset|\geq 3$) that dilates (weakly) $A\in\asetpa$ ($A$ non-extreme) and such that $\lpr(B)>0$, $\forall B\in\bset$.
	Further, suppose that the following $(a)$ and $(b)$ both hold:
	\begin{enumerate}
		\item[(a)] 
		one of the following applies
			\begin{enumerate}
			\item[(a1)]
				$\lpr(A)=0$;
			\item[(a2)]
				%$\forall B\in\bset$ such that $\lpr(A\land B)=0$, $b\leq\frac{1}{P_0(\nega{B})}$;
				$\forall B\in\bset$,
				$\lpr(A\land B)>0,\upr(\nega{A}\land B)>0$;
			\end{enumerate}
		\item[(b)]
		one of the following applies
		\begin{enumerate}
			\item[(b1)]
				$\upr(A)=1$;
			\item[(b2)]
				$\forall B\in\bset$,
				$\lpr(\nega{A}\land B)>0,\upr(A\land B)>0$.
			\end{enumerate}
	\end{enumerate}
Then, there is a (non-trivial) partition other than and coarser than $\bset$, that dilates (weakly) $A$.
\end{theorem}
\begin{proof}
	Let us define
	\begin{align*}
		&\aset^+=\{E\in\asetpa: \lset(A,E)\geq 0\}, \aset^-=\{E\in\asetpa: \lset(A,E)< 0\},\\
		&\bset^+=\{B\in\bset: B\in\aset^+\}, \bset^-=\{B\in\bset: B\in\aset^-\}.
	\end{align*}
	We also write, for notational simplicity,
	\begin{equation*}
	  \bset^+=\{B_1,\ldots,B_k\}, \bset^-=\{B_{k+1},\ldots,B_n\}.
	\end{equation*}
	The proof now continues assuming by contradiction that the thesis is false, i.e. that no partition coarser than $\bset$ dilates $A$.
	
	We investigate three alternatives:
	\begin{enumerate}
		\item[1)]
		$\bset^+\neq\varnothing$, $\bset^-\neq\varnothing$.
		
		Taking $B_i\in\bset^+$, $B_j\in\bset^-$, and defining
		\begin{equation*}
			E_{ij}=\nega{(B_i\lor B_j)}
		\end{equation*}
	we shall reach a contradiction proving that $E_{ij}\notin\aset^+$ and that $E_{ij}\notin\aset^-$.
	\begin{enumerate}
		\item[a)]
		Suppose first that $E_{ij}\in\aset^+$. We consider two distinct partitions coarser than $\bset$.
		\begin{enumerate}
		\item[\ ]
		\emph{Step 1}
		
		Define the partition
		\begin{equation*}
			\bset^\prime=\bset\setminus\{B_i,B_j\}\cup\{B_i\lor B_j\},
		\end{equation*}
	which differs from $\bset$ because it groups together $B_i$ and $B_j$.
	
	$\bset^\prime$ is coarser than $\bset$, and as such does not dilate $A$;
	yet, any of its atoms satisfies Equation \eqref{eq:weakly_dilation}, but for $B_i\lor B_j=\nega{E_{ij}}$.
	
	Preliminarily, note that $\lpr(\nega{E_{ij}})>0$.
	Then, if $(a1)$ holds, we apply Proposition \ref{pro:elle_use} $(a1)$ (to $\bset^\prime$) to get $\lpe(A|\nega{E_{ij}})\leq\lpr(A)$ (since $E_{ij}\in\aset^+$).
	If $(a2)$ holds,
	$\lpr(A\land \nega{E_{ij}}),\upr(\nega{A}\land \nega{E_{ij}})>0$ and we apply now Proposition \ref{pro:elle_use} $(a2)$, getting again $\lpe(A|\nega{E_{ij}})\leq\lpr(A)$.
	Hence, necessarily, $\upe(A|\nega{E_{ij}})<\upr(A)$, because $\nega{E_{ij}}$ does not satisfy \eqref{eq:weakly_dilation}.
	Therefore, straightforwardly if $(b1)$ applies, while, if $(b2)$ holds, recalling that $\lpr(\nega{A}\land \nega{E_{ij}}),\upr(A\land \nega{E_{ij}})>0$, we can conclude by Proposition \ref{pro:elle_use} $(b)$ that
	\begin{equation}
		\label{eq:Eij_neg}
		\lset(\nega{A},E_{ij})<0.
	\end{equation}
	\end{enumerate}
	\begin{enumerate}
	\item[\ ]
	\emph{Step 2}
	
	Now take partition $\bset^{''}=\{B_i,\nega{B_i}\}$.
	Being coarser than $\bset$, $\bset^{''}$ does not dilate $A$.
	This is due to $\nega{B_i}$, since Equation \eqref{eq:weakly_dilation} obtains for $B_i$.
	Similarly to \emph{Step 1}, since $B_i\in\aset^+$, by Proposition~\ref{pro:elle_use} $(a)$ we get $\lpe(A|\nega{B_i})\leq\lpr(A)$, hence $\upe(A|\nega{B_i})<\upr(A)$ necessarily holds.
	Then, Proposition \ref{pro:elle_use} $(b)$ implies
	\begin{equation}
		\label{eq:Bi_neg}
		\lset(\nega{A},B_{i})<0.
	\end{equation}
	\end{enumerate}
	Let us now focus on $B_j$. It satisfies Equation \eqref{eq:weakly_dilation} and this implies by Proposition \ref{pro:elle_use} $(b)$ that $\lset(\nega{A},\nega{B_j})\geq 0$.

	However, since $\nega{B_j}=E_{ij}\lor B_i$, by Proposition \ref{elle_properties}, \eqref{eq:Eij_neg} and \eqref{eq:Bi_neg} we obtain instead
	\begin{equation*}
		\lset(\nega{A},\nega{B_j})=\lset(\nega{A},E_{ij})+\lset(\nega{A},B_i)<0,
	\end{equation*}
	a contradiction. Thus $E_{ij}\notin\aset^+$.
	\item[b)]
	Suppose then that $E_{ij}\in\aset^-$.
	
	Because $B_i$ satisfies Equation \eqref{eq:weakly_dilation}, from Proposition \ref{pro:elle_use} $(a)$ we have that $\lset(A,\nega{B_i})\geq 0$.
	
	However, since $\nega{B_i}=E_{ij}\lor B_j$, from Proposition \ref{elle_properties} we obtain
	\begin{equation*}
		\lset(A,\nega{B_i})=\lset(A,E_{ij})+\lset(A,B_j)<0,
	\end{equation*}
	a contradiction ($\lset(A,E_{ij})<0$ because $E_{ij}\in\aset^-$, $\lset(A,B_j)<0$ because $B_j\in\bset^-$).
	Hence, $E_{ij}\notin\aset^-$.	
	\end{enumerate}
	\item[2)]
	$\bset^+=\varnothing$, $\bset^-=\bset$.
	
	Then, $\lset(A,B_i)<0$, $\forall B_i\in\bset$.
	Apply Corollary \ref{cor:elle_properties} getting
	\begin{equation}
		\label{eq:sumA_nneg}
		0>\sum_{i=1}^n\lset(A,B_i)\geq P_0(A)-\lpr(A)\geq 0,
	\end{equation}
	again a contradiction.
	\item[3)]
	$\bset^+=\bset$, $\bset^-=\varnothing$.
	
	Take any $B\in\bset$. The partition $\{B,\nega{B}\}$, being coarser than $\bset$, does not dilate $A$.
	Hence, we can proceed as in \emph{Step 2} (let $B_i=B$ there) and conclude that
	\begin{equation}
		\label{eq:Aneg_neg}
		\lset(\nega{A},B)<0, \forall B\in\bset.
	\end{equation}
	From \eqref{eq:Aneg_neg} and Corollary \ref{cor:elle_properties}, the following contradiction occurs:
	\begin{equation*}
		%\label{eq:sumnegA_npos}
		0>\sum_{i=1}^n\lset(\nega{A},B_i)\geq P_0(\nega{A})-\lpr(\nega{A})\geq 0.
	\end{equation*}
	\end{enumerate}
\end{proof}
\begin{remark}
\emph{Comments on Theorem \ref{thm:coarsening_basic}.}

If partition $\bset$ dilates strictly $A$, $(a1)$ and $(b1)$ in the statement of Theorem~\ref{thm:coarsening_basic} never apply.
Keeping its remaining assumptions, the proof of Theorem \ref{thm:coarsening_basic} can be adapted to show that  there exists a partition coarser of $\bset$ that dilates strictly $A$.

In both cases Theorem \ref{thm:coarsening_basic} gives a sufficient condition for coarsening.
The condition is not necessary, cf. the later Example \ref{ex:coarser_partition}.
Yet, Theorem  \ref{thm:coarsening_basic} offers some insight on what situations typically ensure the coarsening property.
To see this, suppose firstly that $0<\lpr(A)\leq\upr(A)<1$.
Then, strict positivity of $\lpr(A\land B)$, $\lpr(\nega{A}\land B)$, $\forall B\in\bset$ is enough for coarsening.
In terms of the VBM parameters $a$, $b$, this is equivalent to (recall \eqref{eq:Nearly_Linear_cofficients}, \eqref{eq:Vertical_Barrier_Model_1})
\begin{equation}
\label{eq:param_coars}
(0\leq)-\frac{a}{b}<\min_{B\in\bset}\{\min\{P_0(A\land B),P_0(\nega{A}\land B)\}\}.
\end{equation}
Hence, coarsening arises if the ratio $-\frac{a}{b}$ is `sufficiently low'. Allowing $\lpr(A)=0$ or $\upr(A)=1$ relaxes the constraints in \eqref{eq:param_coars}.
For instance, if $\lpr(A)=0$ Theorem \ref{thm:coarsening_basic} applies when $\lpr(\nega{A}\land B)>0$, $\upr(A\land B)>0$ for all $B$.
However, if further $c\neq 0$ (i.e. $a+b<1$), then it is guaranteed that $\upr(A\land B)>0$ (see also Figure \ref{fig:HBM_VBM_behaviour}, $1)$).
Hence, it only remains to check that $\lpr(\nega{A}\land B)>0$ and consequently \eqref{eq:param_coars} boils down to
\begin{equation*}
-\frac{a}{b}<\min_{B\in\bset}\{P_0(\nega{A}\land B)\}.
\end{equation*}
\end{remark}
Another question is how many partitions coarser than $\bset$, the partition dilating $A$, still dilate $A$.
In the hypotheses of Theorem \ref{thm:coarsening_basic} there is at least one.
The next proposition gives a sufficient condition ensuring dilation with any coarser partition.
\begin{proposition}
	\label{pro:always_dilation}
	Given a VBM $(\lpr,\upr)$ on $\asetpa$, let $\bset\subset\asetpa\setminus\{\emptyset\}$ ($|\bset|\geq 3$) be a finite partition. Suppose that, $\forall B\in\bset$,
	\begin{align}
		\label{eq:non_corr_condition}
		&\lpr(A\land B), \lpr(\nega{A}\land B)>0\text{ and }\\
		\label{eq:non_correlation}
		&P_0(A\land B)=P_0(A)\cdot P_0(B).
	\end{align}
	Then,
	\begin{itemize}
		\item[(i)] $\bset$ weakly dilates $A$;
		\item[(ii)] every partition $\bset^\prime$ coarser than $\bset$ weakly dilates $A$.
	\end{itemize}
\end{proposition}
\begin{proof}
\begin{itemize}
	\item[(i)]
	Equations \eqref{eq:non_corr_condition}, \eqref{eq:non_correlation} ensure that Proposition \ref{pro:dilation_NLM} $2)$ applies. In particular, Equation \eqref{eq:non_correlation} implies that $A$ is $P_0$-non-correlated with any $\nega{B}$ (cf. Remark \ref{rem:dil_suff} $c)$). 
	\item[(ii)]
	Take a generic $B^\prime\in\bset^\prime$.
	Since obviously there is $B\in\bset$ such that $B\Rightarrow B^\prime$, conditions \eqref{eq:non_corr_condition} obtain also when $B$ is replaced by $B^\prime$.
	Further,
	\begin{align*}
		P_0(A\land B')&=P_0(A\land (\vee_{B\in\bset, B\Rightarrow B'} B))=\sum_{B\in\bset, B\Rightarrow B'}P_0(A\land B)\\
		&=\sum_{B\in\bset, B\Rightarrow B'}P_0(A)P_0(B)=P_0(A)P_0(B').
	\end{align*}
	We may therefore apply
	%Proposition \ref{pro:dilation_NLM}
	$(i)$ to $\bset^\prime$, concluding that $\bset^\prime$ dilates $A$.
\end{itemize}	
\end{proof}
The derivation in the proof of Proposition \ref{pro:always_dilation} exploits a simple propagation property of non-correlation, which is reminiscent of similar properties for stochastic independence (see e.g. \cite[Section 4.2]{vicig_epistemic_2000}).
Together with Remark \ref{rem:dil_suff} $c)$, Proposition \ref{pro:always_dilation} highlights the strong role of $P_0$-non-correlation in favouring dilation:
it also prevents getting rid of dilation through any coarsening of the initial partition.

In general, however, $P_0$-non-correlation is not necessary for any partition coarser than $\bset$ to also dilate $A$, see the following example.
\begin{example}
	\label{ex:coarser_partition}
	Continuing Example \ref{ex:dilation}, it may be seen that each of the three binary partitions coarser than $\bset=\{B_1, B_2, B_3\}$ dilates $A$.
	
	Taking for instance $\bset_3=\{B_1\lor B_2, B_3\}$, to verify Equation \eqref{eq:weakly_dilation} we only have to check $B_1\lor B_2$ ($B_3$ already was in Example \ref{ex:dilation} using Proposition~\ref{pro:dilation_NLM}).
	For this we can apply Proposition~\ref{pro:dilation_NLM} again or alternatively \eqref{eq:2_monotone_natural_extension_inf}, \eqref{eq:2_monotone_natural_extension_sup}.
	Following this latter way, we obtain
		\begin{align*}
			\lpe(A|B_1\lor B_2)&=\frac{\lpr(\omega_2\lor\omega_4)}{\lpr(\omega_2\lor\omega_4)+\upr(\omega_1\lor\omega_3)}=\frac{13}{45}=0.2\overline{8}<0.405=\lpr(A);\\
			\upe(A|B_1\lor B_2)&=\frac{\upr(\omega_2\lor\omega_4)}{\upr(\omega_2\lor\omega_4)+\lpr(\omega_1\lor\omega_3)}=\frac{43}{45}=0.9\overline{5}>0.705=\upr(A).
		\end{align*}
	Thus, $\bset_3$ dilates $A$.
	Analogous computations show that $A$ is dilated also by the other two binary partitions coarser than $\bset$.
	Note anyway that Proposition~\ref{pro:always_dilation} does not apply, nor does (since $\lpr(A\land B_2)=0<\lpr(A)$) Theorem \ref{thm:coarsening_basic}.
\end{example}

\subsection{The Extent of Dilation}
\label{subsec:extent_dilation}
Given two events $A$, $B$,
the \emph{imprecision variation} $\Delta(A,B)$ of the uncertainty evaluation for $A$ is defined as
\begin{equation}
	\label{eq:imprec_var}
	\Delta(A,B)=(\upe(A|B)-\lpe(A|B))-(\upr(A)-\lpr(A)).
\end{equation}
The number $\Delta(A,B)$ has a straightforward meaning, and is instrumental in defining the extent of dilation.
When $\Delta(A,B)>0$, it is also called \emph{imprecision increase}.
\begin{definition}
	\label{def:ext_dil}
	Given an event $A$ and a partition $\bset$, the \emph{extent of dilation} $\Delta(A,\bset)$ is defined by
	\begin{align}
		\label{eq:extent_dilation}
		\begin{split}
			\Delta(A,\bset)&=\inf_{B\in\bset}\{\upe(A|B)-\lpe(A|B)-(\upr(A)-\lpr(A))\}\\
			&=\lpr(A)-\upr(A)+\inf_{B\in\bset}\{\upe(A|B)-\lpe(A|B)\}.
		\end{split}
	\end{align}
\end{definition}
Definition \ref{def:ext_dil} corresponds to that in \cite{herron_divisive_1997}, except that it is not restricted to finite partitions.
Results for the extent of dilation were achieved in \cite{herron_divisive_1997} for some special VBMs (the $\varepsilon$-contamination model and the Total Variation Model, both in a finite setting).

We investigate here the extent of dilation for the VBM.
To begin with, the next considerations are helpful in fixing precisely the hypotheses to be assumed in this section.
\begin{enumerate}
	\item[(i)]
	The measure $\Delta(A,\bset)$ may not evaluate adequately weak dilation, since it may be equal to $0$ in such a case, while it may be negative in absence of dilation.
	Therefore, we shall suppose that $\bset$ \emph{strictly dilates} $A$.\footnote
	{
	We are not aware of other measures explicitly devised for graduating the extent of dilation.
	Let us just note in passing that when dilation does not occur, the measure $\Delta^*(A,\bset)=\inf_{B\in\bset}\{\max\{\Delta(A,B),0\}\}$ informs us about the \emph{extent of imprecision increase}, while $\Delta^*(A,\bset)=\Delta(A,\bset)$ under dilation.
	}
	\item[(ii)]
	%When partition $\bset$ is infinite, the infimum in Equation \eqref{eq:extent_dilation} is achieved for any VBM other %than the $\varepsilon$-contamination model.
	
	%To see this, note first that 
	When $\lpr(B)=0$ we have $\Delta(A,B)=1-\upr(A)+\lpr(A)=\max_{B\in\bset}\Delta(A,B)$ (apply Proposition~\ref{pro:natural_extension_case_zero_event} and assumption $(A1)$).
	
	Thus, the infimum in \eqref{eq:extent_dilation} is determined by the remaining atoms of $\bset$, those in
	\begin{equation*}
		\bset_{>0}=\{B\in\bset:\lpr(B)>0\},
	\end{equation*}
if $\bset_{>0}\neq\varnothing$, is trivially equal to $1-\upr(A)+\lpr(A)$ otherwise.
\end{enumerate}
To determine $\Delta(A,\bset)$, we shall consider more subcases, according to whether $\lpe(A|B)=0$, $\upe(A|B)=1$, or both $\lpe(A|B), \upe(A|B)\in ]0,1[$.
Define for this
\begin{align}
	\bset_A^+&=\{B\in\bset_{>0}: \lpe(A|B),\upe(A|B)\in ]0,1[\},\\
	\bset_A^1&=\{B\in\bset_{>0}: \lpe(A|B)\in ]0,1[,\upe(A|B)=1\},\\
	\bset_A^0&=\{B\in\bset_{>0}: \lpe(A|B)=0,\upe(A|B)\in ]0,1[\}.
\end{align}
The extent of dilation is then given by the following
\begin{theorem}
	\label{thm:extent_dilation_VBM}
	Let $(\lpr,\upr)$ be a VBM on $\asetpa$.
	Let $\bset\subset\asetpa\setminus\{\emptyset\}$ be a partition that strictly dilates $A\in\asetpa$.
	Then, the extent of dilation is
	\begin{equation}
		\label{eq:extent_dilation_VBM}
		\Delta(A,\bset)=\lpr(A)-\upr(A)+\min\left\{1,\upe(A|B^*)-\lpe(A|B^*),\frac{1}{1+M_0},\frac{1}{1+M_1}\right\},
	\end{equation}
	where $B^*$ is such that
	\begin{equation}
		\label{eq:B*}
		P_0(B^*)=\max_{B\in\bset_A^+} P_0(B)
	\end{equation}
	and
	\begin{equation*}
		\label{eq:M0M1}
		M_0=\sup_{B\in\bset_A^0}\frac{\lpr(\nega{A}\land B)}{\upr(A\land B)},\ M_1=\sup_{B\in\bset_A^1}\frac{\lpr(A\land B)}{\upr(\nega{A}\land B)}.
	\end{equation*}
It is understood that if any among $\bset_A^+$, $\bset_A^0$, $\bset_A^1$ is empty the corresponding term in Equation \eqref{eq:extent_dilation_VBM} has to be skipped. 
\end{theorem}
\begin{proof}
	The term $1-\upr(A)+\lpr(A)$, obtained when the minimum in \eqref{eq:extent_dilation_VBM} is $1$,  covers the case that $\bset_{>0}=\varnothing$, according to $(ii)$, and that
	$\bset_{>0}\neq\varnothing$ but $\lpe(A|B)=0$, $\upe(A|B)=1$, $\forall B\in\bset$. 
	Note that when $B\in\bset_{>0}$, if $\lpe(A|B)=0$ then $\upe(A|B)>0$ ($\upe(A|B)=0$ would prevent strict dilation),
	and similarly $\upe(A|B)=1$ implies $\lpe(A|B)<1$.
	
	In the remaining alternatives, we determine the infimum of $\Delta(A,B)$ for $B\in\bset_A^+$, $B\in\bset_A^1$, $B\in\bset_A^0$ as follows.
	\begin{itemize}
		\item[(a)] 
		Let $B\in\bset_A^+$.
		
		From \eqref{eq:NE_VBM_inf}, \eqref{eq:NE_VBM_sup}, \eqref{NE_VBM_coefficients}
		\begin{align}
			\nonumber
			\upe(A|B)-\lpe(A|B)&=b_B P_0(A|B)+c_B -b_B P_0(A|B)-a_B\\
			\label{eq:2a2b}
			&=1-2a_B-b_B\\
			\nonumber
			&=1-\frac{2a+b P_0(B)}{b P_0(B)+1-b}.
		\end{align}
	To minimise $\upe(A|B)-\lpe(A|B)$, maximise $\phi(x)=\frac{2a+bx}{bx+1-b}$.
	We have that $\phi^\prime(x)=\frac{b(1-b-2a)}{(bx+1-b)^2}\geq 0$,
	because $b>0$, $b+2a\leq b+a\leq 1$ in the VBM. Anyway, $\phi^\prime(x)=0$ iff  $b+2a=1$ iff $a=c$, i.e. iff $\lpr=\upr$ in the VBM, a limit situation that cannot originate dilation and may therefore be discarded.
	This means that $\phi$ is strictly increasing, and such is its restriction on $\{P_0(B):B\in\bset_A^+$\}.
	Therefore, recalling \eqref{eq:B*},
	\begin{equation*}
		\inf_{B\in\bset_A^+}\{\upe(A|B)-\lpe(A|B)\}=\upe(A|B^*)-\lpe(A|B^*).
	\end{equation*}
	\item[(b)]
	Let $B\in\bset_A^0$.
	
	In this case $\upe(A|B)>0$ ensures by \eqref{eq:2_monotone_natural_extension_sup} that $\upr(A\land B)>0$.
	Using this and \eqref{eq:2_monotone_natural_extension_sup} again,
	\begin{equation*}
		\inf_{B\in\bset_A^0}\{\upe(A|B)-\lpe(A|B)\}=\inf_{B\in\bset_A^0}\frac{\upr(A\land B)}{\upr(A\land B)+\lpr(\nega{A}\land B)}=\frac{1}{1+M_0}.
	\end{equation*}
	\item[(c)]
	Let $B\in\bset_A^1$.
	
	$\lpe(A|B)\in ]0,1[$ guarantees by \eqref{eq:2_monotone_natural_extension_inf} that $\upr(\nega{A}\land B)>0$.
	Because of this, using firstly conjugacy of $\lpe$, $\upe$ and \eqref{eq:2_monotone_natural_extension_sup}, we write
	\begin{align*}
		\inf_{B\in\bset_A^1}\{\upe(A|B)-\lpe(A|B)\}&=\inf_{B\in\bset_A^1}\{1-\lpe(A|B)\}
		=\inf_{B\in\bset_A^1}\{\upe(\nega{A}|B)\}\\
		&=\inf_{B\in\bset_A^1}\frac{\upr(\nega{A}\land B)}{\upr(\nega{A}\land B)+\lpr(A\land B)}=\frac{1}{1+M_1}.
	\end{align*}
	\end{itemize}
Equation \eqref{eq:extent_dilation_VBM} then follows from the previous derivations.
\end{proof}

When the VBM is not an $\varepsilon$-contamination model, i.e. when $a\neq 0$ in Definition \ref{def:Vertical_Barrier}, the infimum in Equation \eqref{eq:extent_dilation} is achieved, since then $\bset_{>0}$ (when non-empty) is finite, and so are its subsets $\bset_{A}^{+}$, $\bset_{A}^{1}$, $\bset_{A}^{0}$.
In fact, from Definition \ref{def:Vertical_Barrier}, $\lpr(B)>0$ iff $P_0(B)>-\frac{a}{b}>0$, a condition that may obtain for finitely many $B\in\bset$. Consequently, every infimum in the proof of Theorem~\ref{thm:extent_dilation_VBM} is a minimum.

We point out that when the relevant evaluations are non-extreme, i.e. when $\lpr(B)>0$ and $\lpe(A|B), \upe(A|B)\in ]0,1[, \forall B\in\bset$, the computation of $\Delta(A,\bset)$ simplifies to $\upe(A|B^*)-\lpe(A|B^*)$,
with $B^*$ easily detected by means of Equation \eqref{eq:B*}.

Still referring to this non-extreme situation, in general dilation (weak or strict) may or may not occur.
It is anyway interesting to observe that the weaker phenomenon of \emph{imprecision increase} always obtains when $b<1$ in the VBM.
This is stated precisely in the next
\begin{proposition}
	\label{pro:impr_incr_bmin1}
	Given a VBM on $\asetpa$, $A\in\asetpa$ such that $\lpr(A),\upr(A)\in ]0,1[$
	and a partition $\bset\subset\asetpa\setminus\{\emptyset\}$, suppose that, $\forall B\in\bset$, $\lpr(B)>0$, $P_0 (B)<1$, $\lpe(A|B), \upe(A|B)\in ]0,1[$. Let $b<1$ in the VBM.
	Then,
	\begin{equation*}
		\Delta(A,B)>0, \forall B\in\bset.
	\end{equation*}
\end{proposition}
\begin{proof}
	Recalling \eqref{eq:Nearly_Linear_cofficients}, \eqref{eq:Vertical_Barrier_Model_1}, \eqref{eq:Vertical_Barrier_Model_2}, \eqref{eq:imprec_var} and \eqref{eq:2a2b},
	\begin{equation*}
		\Delta(A,B)=bP_0(A)+a-(bP_0(A)+c)+1-2a_B-b_B=2(a-a_B)+b-b_B.
	\end{equation*}
	From \eqref{NE_VBM_coefficients},
	\begin{enumerate}
		\item[$\bullet$]
		$a-a_B=a(1-\frac{1}{bP_0(B)+1-b})\geq 0$, because $bP_0(B)+1-b=1-bP_0(\nega{B})<1$ ($P_0(\nega{B})>0$), while $a\leq0$;
		\item[$\bullet$]
		$b-b_B=b(1-\frac{P_0(B)}{bP_0(B)+1-b})>0$ iff $\frac{P_0(B)}{bP_0(B)+1-b}<1$ iff $P_0(B)(1-b)<1-b$,
		which holds because $b<1$ and $P_0(B)<1$.
	\end{enumerate}
	Thus, $\Delta(A,B)>0$.
\end{proof}

\subsection{Constriction}
\label{subsec:constriction}
Following \cite{herron_divisive_1997}, we have that
\begin{definition}
\label{def:constriction}
A partition $\bset$ \emph{constricts} $A$ iff
\begin{equation}
	\label{eq:constriction}
	\lpr(A)\leq\lpe(A|B)\leq\upe(A|B)\leq\upr(A), \forall B\in\bset,
\end{equation}
with at least one of the outer inequalities in \eqref{eq:constriction} being strict, for at least one $B\in\bset$.
\end{definition}
Clearly, constriction is a sort of opposite of dilation.
Guaranteeing that the imprecision of the conditional evaluation will never be larger than that of the unconditional one,
and at least once smaller, is a very desirable property.
Unfortunately, it is quite uncommon, being subject to rather restrictive conditions.

Constriction was explored in \cite{herron_divisive_1997} for the $\varepsilon$-contamination model, proving that it never takes place if $\varepsilon>0$.

We shall show hereafter that it occurs for no VBM either, under `normal' situations.
Basically, `normal' means that constriction requires some relevant evaluations to be extreme (a situation generally not considered in \cite{herron_divisive_1997}) in order to take place.
Even this may be far from sufficient:
a thorough investigation of constriction shows that it occurs in very special situations only.
However, we shall not describe in detail all the very specific cases that originate or not constriction,
which might be somewhat tedious to the reader.
The results we report should already give a clear idea of how restrictive constriction is.

Studying when $\bset$ constricts $A$, we shall suppose that
\begin{itemize}
\item[\textbf{(A4)}] 
$\lpr(A)<\upr(A)$.
\end{itemize}
Note that $(A4)$, which is necessary to `give a chance' to constriction, is stronger than $(A2)$.

Let us now turn to the case that $\lpr(B)>0$, $\forall B\in\bset$, i.e. that $\bset_0=\varnothing$, with
\begin{equation*}
	\bset_0=\{B\in\bset:\lpr(B)=0\}.
\end{equation*}
Then, we have that
\begin{proposition}
	\label{pro:constriction}
	Given a VBM $(\lpr,\upr)$ on $\asetpa$, take $A\in\asetpa$ and partition $\bset\subset\asetpa\setminus\{\emptyset\}$.
	Suppose $\bset_0=\varnothing$.
	If there is $B^*\in\bset$ such that $\lpr(A\land B^*), \lpr(\nega{A}\land B^*)>0$, then $\bset$ does not constrict $A$.
\end{proposition}
\begin{proof}
	We prove that \eqref{eq:constriction} does not apply with $B=B^*$.
	
	Preliminarily, we deduce easily from Lemma \ref{lem:partial_dilation} that in the current hypotheses
	\begin{align}
		\label{co_prel_1}
		\lpe(A|B^*)&>\lpr(A) \text{ iff } \lpr(A)>P_0(A|\nega{B^*}),\\
		\label{co_prel_2}
		\lpe(A|B^*)&=\lpr(A) \text{ iff } \lpr(A)=P_0(A|\nega{B^*}),\\
		\label{co_prel_3}
		\upe(A|B^*)&<\upr(A) \text{ iff } \upr(A)<P_0(A|\nega{B^*}),\\
		\label{co_prel_4}
		\upe(A|B^*)&=\upr(A) \text{ iff } \upr(A)=P_0(A|\nega{B^*}).
	\end{align} 
	Then, we see firstly that the outer inequalities in \eqref{eq:constriction} cannot both be equalities, since then $\lpr(A)=\upr(A)$ by \eqref{co_prel_2}, \eqref{co_prel_4}, against assumption $(A4)$.
	Secondly, if at least one of the outer inequalities in \eqref{eq:constriction} is strict,
	recalling also Equations \eqref{co_prel_1}, \eqref{co_prel_3}, we deduce that $\lpr(A)>\upr(A)$. This conflicts with coherence of $\lpr$, $\upr$.
\end{proof}
Proposition \ref{pro:constriction} lets us already see that constriction is not the rule with VBMs.
In fact, when $\lpr(B)>0$, $\forall B\in\bset$, it says that there is no constriction if (recall \eqref{eq:2_monotone_natural_extension_inf}, \eqref{eq:2_monotone_natural_extension_sup}) there is $B^*$ such that $\lpe(A|B^*)$, $\upe(A|B^*)$ are non-extreme, i.e. take values in $]0,1[$. This is anyway a usual situation.

Thus, we are bound to allow extreme evaluations or drop some of the standard assumptions in the search for constriction. The next proposition illustrates this.
\begin{proposition}
	\label{pro:constr_extr}
	Given a VBM $(\lpr,\upr)$ on $\asetpa$, suppose now $\bset_0\neq\varnothing$.
	If there exists $B^*\in\bset_0$ such that $B^*\not\Rightarrow A$, $B^*\not\Rightarrow\nega{A}$, $\bset$ constricts $A$ iff:
	$\lpr(A)=0$, $\upr(A)=1$, and there is $B'\in\bset_0$ such that either $B'\Rightarrow A$ or $B'\Rightarrow\nega{A}$.
\end{proposition}
\begin{proof}
	From Proposition \ref{pro:natural_extension_case_zero_event}, $\lpe(A|B^*)=0$, $\upe(A|B^*)=1$,
	thus $B^*$ satisfies (weakly) Equation \eqref{eq:constriction} iff $\lpr(A)=0$, $\upr(A)=1$.
	
	These two extreme conditions are therefore necessary for constriction.
	Assuming them, we see that $\lpe(A|B)=0$, $\upe(A|B)=1$ for any $B\in\bset\setminus\bset_0$.
	This is because $\lpr(A)=0$, $\upr(A)=1$ imply $\lpr(A\land B)=\lpr(\nega{A}\land B)=0$, and by \eqref{eq:2_monotone_natural_extension_inf}, \eqref{eq:2_monotone_natural_extension_sup} $\lpe(A|B)=0$, $\upe(A|B)=1$.
	Thus, \eqref{eq:constriction} holds weakly for any such $B$.
	
	To find some conditioning event such that at least one of the outer inequalities in \eqref{eq:constriction} is strict,
	it is necessary that for some $B'\in\bset_0$ either $B'\Rightarrow A$ or $B'\Rightarrow\nega{A}$.
	This is also sufficient:
	\begin{itemize}
		\item[-] if $B'\Rightarrow A$, then $\lpe(A|B')=\upe(A|B')=1$ and the leftmost inequality in \eqref{eq:constriction} is strict;
		\item[-] if $B'\Rightarrow \nega{A}$, then $\lpe(A|B')=\upe(A|B')=0$ and the rightmost inequality in \eqref{eq:constriction} is strict.
	\end{itemize}
\end{proof}
We point out that constriction is achieved in Proposition \ref{pro:constr_extr} only introducing extreme evaluations and \emph{dropping} logical independence of $A$ from $\bset$, i.e. dropping assumptions $(A1)$ and $(A3)$ in Section \ref{subsec:dilation_with_NL}.

Concerning logical independence, an immediate follow-up of Proposition~\ref{pro:constr_extr} is
\begin{corollary}
	\label{cor:constriction}
	If, given a VBM, $A$, $\bset$, it holds that $\bset_0\neq\varnothing$ and $A$ is logically independent of $\bset$, then $\bset$ does not constrict $A$.
\end{corollary}
Further conditions that ensure or (more commonly) prevent constriction may be found, still giving up some of the axioms $(A1)\div(A3)$.
We shall not enter details of these very specific situations.

\section{Conclusions}
\label{sec:conclusions}
The present investigation confirms that coherent NL models are a manageable family of imprecise probabilities also with respect to conditioning: simple formulae are available for extending them on $\asetpa|B$.
VBMs have a prominent role, since they are also stable, i.e. return a conditional model of the same kind.

While studying dilation with NL models, we have also undertaken a more general investigation of the hypotheses that 
make it possible to decide whether dilation occurs or not.
It is the case of extreme events $A$, or of various dependence structures relating $A$ with partition $\bset$.
Ruling out these cases we are left with assumptions still more general than often assumed.
Thus, logical independence of $A$ from $\bset$ does not imply the common positivity hypothesis $\lpr(B)>0$, $\forall B\in\bset$.
Nevertheless, dilation for coherent NL models can be characterized.

We have also seen how to apply an existing measure of dilation, the extent of dilation, and that the opposite phenomenon of constriction, although desirable, remains quite rare.
Changing the partition $\bset$ may be a way to prevent dilation. Choosing a more refined partition tends to introduce some logical dependencies and hence to rule out dilation, while increasing the number of evaluations to work with.
As for the opposite choice, we have seen that not every coarser partition is dilation-immune, at least in the assumptions of Theorem \ref{thm:coarsening_basic}.
It is true anyway that this result, although generalising previous ones, is only sufficient for coarsening.
A general characterisation of coarsening is still an open problem, to the best of our knowledge.

In general, little is so far known about ways of avoiding dilation, and also on the weaker phenomenon of imprecision increase.
These could be promising issues for future work.

\section{Appendix}
\label {sec:appendix}

\subsection{Conditioning with the Regular Extension}
\label{subsec:cond_RE}

It appears also from Proposition \ref{pro:natural_extension_case_zero_event} that
inferences with the natural extension may lead to vacuous assessments, when $\lpr(B)=0$.
We discuss here the alternative \emph{regular extension}, suggested in the literature for these cases.

\subsubsection{W-coherence of the Regular Extension}
\label{subsub:W-coh_reg_ext}
We begin in a general framework, not limited to NL models. Let for this $\lpr$ be a coherent lower probability on $\asetpa$, $\upr$ its conjugate, $\mset$ their credal set.\footnote{
	More generally, the results in this section apply when $\lpr$ and $\upr$ are defined on an algebra of events.
}
Take $B\in\asetpa\setminus\{\emptyset\}$, and suppose for the moment that
\begin{eqnarray}
	\label{eq:reg_condition}
	\lpr(B)=0<\upr(B).
\end{eqnarray}
Define further
\begin{eqnarray*}
	\label{eq:reg_credal_sets}
	\mset^+&=&\{P\in\mset\text{ s.t. }P(B)>0\},\\
	%\mset^+_B=\{P:\asetpa\cup\asetpa|B\rightarrow\rset: P\mathord\restriction_{\asetpa}\in\mset^+\}.\\
	\mset^+_B&=&\{P:\asetpa\cup\asetpa|B\rightarrow\rset\text{ s.t. }P|_{\asetpa}\in\mset^+\}.
\end{eqnarray*}
In words, $\mset^+_B$ is given by the extensions of the (unconditional) probabilities in $\mset^+$ to the events $A|B$, with $A$ varying in $\asetpa$.
%For each $P\in\mset^+$, its extension is unique and given by $P(A|B)=\frac{P(A\land B)}{P(B)}$, $\forall A\in\asetpa$.
Note that Equation~\eqref {eq:reg_condition} and the Envelope Theorem guarantee that $\mset^+\neq\varnothing$.

Given this, define the lower probability on $\asetpa\cup\asetpa|B$
\begin{eqnarray}
	\label{eq:reg_inf}
	\lreg=\inf_{P\in\mset^+_B} P.
\end{eqnarray}
\begin{proposition}
	\label{pro:reg_coherent}
	Given a coherent lower probability $\lpr:\asetpa\rightarrow\rset$, such that \eqref{eq:reg_condition} applies,
	\begin{enumerate}
		\item[(a)]
		$\lreg:\asetpa\cup\asetpa|B\rightarrow\rset$ is a W-coherent lower probability.
		\item[(b)] On $\asetpa$, $\lreg=\lpr$.
	\end{enumerate}
\end{proposition}
\begin{proof}
	\begin{enumerate}
		\item[(a)]
		Follows from \eqref{eq:reg_inf} and the Lower Envelope Theorem.
		\item[(b)]
		$\forall A\in\asetpa$, by definition, $\lreg(A)=\inf\{P(A):P\in\mset^+_B\}=\inf\{P(A):P\in\mset^+\}$.
		Two alternatives arise:
		\begin{enumerate}
			\item[(b1)]
			If $\exists P^*\in\mset^+: P^*(A)=\lpr(A)$, then obviously $\lreg(A)=\lpr(A)$.
			\item[(b2)]
			Otherwise, if $P^*(B)=0$ $\forall P^*\in\mset$ such that $P^*(A)=\lpr(A)$, then necessarily $P(A)>\lpr(A)$ for any $P\in\mset^+$.
			Take one such $P$.
			Since $\mset$ is convex, for any $\varepsilon \in\ ]0,1[$, $P_\varepsilon=\varepsilon P+(1-\varepsilon)P^*\in\mset$ and also
			$P_\varepsilon\in\mset^+$, because $P_\varepsilon(B)=\varepsilon P(B)>0$.
			We have that
			\begin{eqnarray*}
				P_\varepsilon(A)=\varepsilon P(A)+(1-\varepsilon)\lpr(A)
				=\varepsilon(P(A)-\lpr(A))+\lpr(A).
			\end{eqnarray*}
			This means that $\lpr(A)\leq\lreg(A)=\inf\{P(A):P\in\mset^+\}\leq\inf\{P_\varepsilon (A):\varepsilon\in\ ]0,1[\}=\lpr(A)$, hence $\lreg(A)=\lpr(A)$.
			
			Note that case $(b2)$ applies to $A=B$ too.
		\end{enumerate}
	\end{enumerate}
\end{proof}
At this point, let us remove assumption \eqref{eq:reg_condition} and extend the definition of $\lreg$ in \eqref{eq:reg_inf} as follows:
\begin{definition}
	Given a coherent lower probability $\lpr$ on $\asetpa$,
	its \emph{regular extension} $\lreg$ on $\asetpa\cup\asetpa|B$ is given by
	\begin{align}
		\label{eq:RegExt}
		\lreg=\begin{dcases}
			\inf_{P\in\mathcal{M}^+_B} P  &\text{ if } \mset^+\neq\varnothing\\
			\lpe &\text{ if } \mathcal{M}^+=\varnothing.
		\end{dcases}
	\end{align}
\end{definition}
Note that $\lreg=\lpe$ also when ($\mset^+\neq\varnothing$ and) $\lpr(B)>0$, since then $\mset=\mset^+$.
In fact, this equality implies, $\forall A\in\asetpa$, $\lreg(A)=\lpe(A)$ and, recalling also \eqref{eq:nat_ext_expr},
$\lpe(A|B)=\inf\{P(A|B):P\in\mset^+_B\}=\lreg(A|B)$.

Clearly then, by the properties of the natural extension and Proposition~\ref{pro:reg_coherent},
\emph{the regular extension is W-coherent and $\lreg=\lpe=\lpr$ on $\asetpa$}.

\subsubsection{The Regular Extension of $2$-monotone Models}
\label{subsub:reg_ext_2_mon_mod}
Of course, what primarily matters is to make clear when the regular and the natural extensions differ. Unfortunately, this instance is very limited if $\lpr$ is $2$-monotone:
\begin{proposition}
	\label{pro:lreg_2_mon}
	When $\lpr$ is coherent and $2$-monotone on $\asetpa$, for any $A\in\asetpa$ we have that $\lreg(A|B)\neq\lpe(A|B)$ iff
	\begin{eqnarray}
		\label{2-mon_reg_nat_diff}
		0=\lpr(B)<\upr(B), \nega{A}\land B\neq\emptyset\text{ and }\upr(\nega{A}\land B)=0.
	\end{eqnarray}
	In this case, $\lreg(A|B)=1>0=\lpe(A|B)$.
\end{proposition}
\begin{proof}
	In general, even when $\lpr$ is not $2$-monotone, it is necessary for $\lreg$ to differ from $\lpe$ that $\mset^+\neq\varnothing$ and $\lpr(B)=0$, or equivalently that $\lpr(B)=0<\upr(B)$.
	%However, if this condition occurs but $\lpr(B)>0$, we still have that $\lreg=\lpe$ (since $\mset=\mset^+$, %$\lpe(A|B)=\inf\{P(A|B): P\in\mset^+_B\}$).
	%
	%We may therefore restrict our attention to the case $\mset^+\neq\varnothing$, $\lpr(B)=0$, or equivalently %$\lpr(B)=0<\upr(B)$.
	With $2$-monotonicity the following result, proven in \cite[Lemma 7.1]{walley_coherent_1981}, is helpful:
	
	If $\lpr$ is $2$-monotone on $\asetpa$ and $\lpr(B)=0<\upr(B)$, then, $\forall A\in\asetpa$,
	\begin{align}
		\label{eq:RegExt_Walley}
		\lreg(A|B)=\begin{dcases}
			1 &\text{ if } \upr(\nega{A}\land B)=0\\
			0 &\text{ if } \upr(\nega{A}\land B)>0.
		\end{dcases}
	\end{align}
	Recall also from Proposition \ref{pro:natural_extension_case_zero_event} that, when $\lpr(B)=0$, $\lpe(A|B)=0$ if $B\not\Rightarrow A$ (equivalent to $\nega{A}\land B\neq\emptyset$), while $\lpe(A|B)=1$ when $B\Rightarrow A$, i.e. when $\nega{A}\land B=\emptyset$.
	The thesis then follows as a joint consequence of Proposition~\ref{pro:natural_extension_case_zero_event} and \eqref{eq:RegExt_Walley}.
\end{proof}
The result in Proposition \ref{pro:lreg_2_mon} (which, in the case that $\prt$ is finite, can be deduced also from results in \cite{miranda_coherent_2015}) applies by Proposition \ref{pro:NL_properties} in particular to coherent NL models: the regular extension differs for them from the natural extension only in the very special condition \eqref{2-mon_reg_nat_diff}.
Moreover, in this case $\lreg$ takes the extreme value $1$.

\subsubsection{The Case of Vertical Barrier Models}
\label{subsub:case_VBM}
We can further appreciate how restrictive condition \eqref{2-mon_reg_nat_diff} is for NL models if we focus on VBMs.
In fact, within VBMs condition \eqref{2-mon_reg_nat_diff} applies \emph{only} to PMMs:
\begin{proposition}
	\label{pro:VBM_becomes_PMM}
	A VBM that satisfies condition \eqref{2-mon_reg_nat_diff} is a PMM.
\end{proposition}
\begin{proof}
Conditions $\nega{A}\land B\neq\emptyset$, $\upr(\nega{A}\land B)=0$ in \eqref{2-mon_reg_nat_diff} imply by \eqref{eq:Vertical_Barrier_Model_2} $bP_0(\nega{A}\land B)+c=0$ and hence (since $b>0$ by \eqref{eq:Nearly_Linear_cofficients}) $c=0$.

When $c=0$, condition $\upr(B)>0$, requiring that $bP_0(B)+c=bP_0(B)>0$, guarantees that $P_0(B)>0$.

Finally, when $P_0(B)>0$ condition $\lpr(B)=0$, equivalent to $P_0(B)\leq -\frac{a}{b}$, ensures that $a<0$, while (by \eqref{eq:Nearly_Linear_cofficients}) $b=1-a>1$.

Conditions $c=0$, $b>1$, $a<0$ identify a PMM (cf. Section \ref{subsub:Nearly_Linear_Models}).
\end{proof}
Thus, the regular and natural extensions are identical for any VBM that is not a PMM.
%Note also that a VBM may `introduce zeros': there may exist events $B$ such that $\lpr(B)=0$ while $P_0(B)>0$.
%Clearly, here the regular extension would be most useful. Yet, the only VBM that may differentiate the regular extension, %i.e. the PMM, does not `introduce zeros'.

\subsection{Dilation, Logical Independence and Extreme Events}
\label{sec:dilation_generalities}
In this section we discuss the conditions that can be reasonably required in studying dilation.
Reasonably means that if either condition does not apply, we already know what happens to dilation and no further investigation is needed.

We explore a very general situation, not asking that a NL model is assessed but considering just an unconditional coherent assessment which includes $\lpr(A)$, $\upr(A)$, $\lpr(B)$ for given non-trivial events $A$, $B$.

Moreover, dilation may affect any coherent extension of the given assessment to $A|B$, $\forall B\in\bset$---not solely the natural extension considered in Definition \ref{def:weakly_dilation}.
With one such coherent extension, $\lpr(\cdot|\cdot)$, $\upr(\cdot|\cdot)$, weak and strict dilation are defined as in Definition \ref{def:weakly_dilation}, replacing equation \eqref{eq:weakly_dilation} with
\begin{equation}
	\label{eq:dilation_general}
	\lpr(A|B)\leq\lpr(A)\leq\upr(A)\leq\upr(A|B), \forall B\in\bset.
\end{equation}
\subsubsection{Dilation and Logical (In)dependence}
\label{subsubsec:dil_log_ind}
A first relevant issue concerning dilation regards the logical relationships between $A$ and $\bset$.
Precisely,
\begin{definition}
	\label{def:depend}
	Given an event $A$ and a partition $\bset$,
	\begin{enumerate}
		\item[(a)] $A$ is \emph{logically dependent} on $\bset$ if either $B\Rightarrow A$ or $B\Rightarrow\nega{A}$, $\forall B\in\bset$;
		\item[(b)] $A$ is \emph{(logically) semidependent} on $\bset$ if there is $B'\in\bset$ such that  $B'\not\Rightarrow A$ and $B'\not\Rightarrow\nega{A}$, and there is $B''\in\bset$ such that either $B''\Rightarrow A$ or $B''\Rightarrow\nega{A}$.
	\end{enumerate}
\end{definition}
When $A$ is semidependent on $\bset$ it may happen that
\begin{enumerate}
	\item[(b1)] there is no $B\in\bset$ such that $B\Rightarrow A$;
	\item[(b2)] there is no $B\in\bset$ such that $B\Rightarrow \nega{A}$.
\end{enumerate}
In case $(b1)$ (in case $(b2)$) we say that the semidependence of $A$ is \emph{one-sided of type $1$} (\emph{one-sided of type $2$}).

It holds that
\begin{proposition}
	\label{pro:indep_dilation}
	Partition $\bset$ does not dilate weakly (and not trivially) $A$ if $A$ is either logically dependent on $\bset$ or semidependent on $\bset$ and its semidependence is not one-sided. 
	%not logically independent of $\bset$,
	%except for the case that $A$ is semidependent of types $1$ or, alternatively, $2$ on $\bset$ and either %$\lpr(A)=\upr(A)=0$ or $\lpr(A)=\upr(A)=1$, respectively.
\end{proposition}
\begin{proof}
	If $A$ is not logically independent of $\bset$, there is $B\in\bset$ such that either $B\Rightarrow A$ or $B\Rightarrow\nega{A}$.
	\begin{enumerate}
		\item[(i)] When $B\Rightarrow\nega{A}$, then $A\land B=\emptyset$, and consequently coherence requires that $\lpr(A|B)=\upr(A|B)=0$.
		
		Thus $B$ satisfies \eqref{eq:dilation_general} (with equality everywhere) only if $\lpr(A)=\upr(A)=0$.
		\item[(ii)] When $B\Rightarrow A$, then $A|B=B|B$ and $\lpr(A|B)=\upr(A|B)=1$.
		Thus $B$ satisfies \eqref{eq:dilation_general} (with equality everywhere) only if $\lpr(A)=\upr(A)=1$.
	\end{enumerate}
	
	From $(i)$ and $(ii)$ above, we see that weak dilation cannot occur if $A$ is semidependent on $\bset$ and its dependence is not one-sided,
	since $A$ should satisfy the conflicting requirements $\lpr(A)=1$ and $\upr(A)=0$. 
	
	For the same reason logical dependence of $A$ on $\bset$ prevents dilation. In fact, logical dependence should be one-sided too,
	which happens only when $A=\emptyset$ (if $B\Rightarrow\nega{A}$, $\forall B\in\bset$) or $A=\Omega$ (if $B\Rightarrow A$, $\forall B\in\bset$).
\end{proof}

Proposition \ref{pro:indep_dilation} does not treat the case that $A$ is logically semidependent on $\bset$ and its semidependence is one-sided. This is done in the next corollary.

\begin{corollary}
	\label{cor:1-semidependence}
	\begin{enumerate}
		\item[(a)] If $A$ is \emph{semidependent of type $1$} on $\bset$, $\bset$ dilates weakly, but not trivially, $A$ iff:
		\begin{enumerate}
			\item[(a1)] $\lpr(A)=\upr(A)=0$;
			\item[(a2)] $\lpr(A|B)=0, \forall B\in\bset$;
			\item[(a3)] $\exists B\in\bset: \upr(A|B)>0$.
		\end{enumerate}
		\item[(b)] If $A$ is \emph{semidependent of type $2$} on $\bset$, $\bset$ dilates weakly, but not trivially, $A$ iff:
		\begin{enumerate}
			\item[(b1)] $\lpr(A)=\upr(A)=1$;
			\item[(b2)] $\upr(A|B)=1, \forall B\in\bset$;
			\item[(b3)] $\exists B\in\bset: \lpr(A|B)<1$.
		\end{enumerate}
	\end{enumerate}
\end{corollary}
\begin{proof}
	We prove $(a)$, $(b)$ being analogous.
	If $A$ is semidependent of type $1$ on $\bset$, by definition there is $B\in\bset$ such that $B\Rightarrow\nega{A}$.
	From $(i)$ in the proof of Proposition \ref{pro:indep_dilation},
	any such $B$ satisfies \eqref{eq:dilation_general}, trivially, only if $\lpr(A)=\upr(A)=0$.
	In order for weak dilation to hold, it is then necessary that $\lpr(A|B)=0$ also for the remaining atoms $B$ of $\bset$, those not implying $\nega{A}$.
	At this point, weak dilation is not trivial iff $(a3)$ applies.
\end{proof}
Thus, there is uncertainty whether dilation occurs or not only if $A$ is logically independent of partition $\bset$.
\subsubsection{Dilation with the Natural Extension and Extreme Events}
\label{subsubsec:dilation_nat_ext_extr_ev}
Even requiring logical independence of $A$ from $\bset$, a second relevant case occurs when $A$ is an extreme event. Here, we can establish \emph{a priori} that dilation occurs (trivially or not),
but unlike the previous subsection we have to assume that the natural extension is performed.
The result is stated precisely in the next proposition.
Its part $(a)$ generalises \cite[Lemma 3]{pelessoni_dilation_2020}, which applies to NL models only.
\begin{proposition}
	\label{pro:dilation_extreme}
	Given a(n unconditional) coherent lower probability $\lpr$ on $\dset\supset\{A,\nega{A},B\}$ and its conjugate $\upr$, let $A$ be an extreme event logical independent from $\bset$. Apply the natural extension to compute $\lpe(A|B)$, $\upe(A|B)$, $\forall B\in\bset$.
	\begin{enumerate}
		\item[(a)] If $\lpr(B)>0$,  weak dilation occurs trivially (all inequalities are equalities in equation \eqref{eq:weakly_dilation}). 
		\item[(b)] If $\lpr(B)=0$, \eqref{eq:weakly_dilation} applies (not trivially, one outer inequality being strict).
	\end{enumerate}
\end{proposition}
\begin{proof}
As for $(a)$, if $\lpr(A)=\upr(A)=0$, then also $P(A)=0$, for any probability $P$ in the credal set $\mset$.
Hence $P(A|B)=\frac{P(A\land B)}{P(B)}=0$ (since $P(A\land B)=0$ by monotonicity of $P$ and $P(B)\geq\lpr(B)>0$).
It follows, recalling the Envelope Theorem, $\lpe(A|B)=\inf_{P\in\mset} P(A|B)=0$, $\upe(A|B)=\sup_{P\in{\mset}} P(A|B)=0$.
Thus, \eqref{eq:weakly_dilation} holds trivially.
		
If instead $\lpr(A)=\upr(A)=1$, then $\lpe(A|B)=1-\upe(\nega{A}|B)=1$ (the second equality being justified by conjugacy and the argument above). Hence, $\upe(A|B)=1$ and \eqref{eq:weakly_dilation} holds trivially again.

$(b)$ was proven in \cite[Lemma 2]{pelessoni_dilation_2020}.
\end{proof}

%\paragraph{\textbf{\emph{Acknowledgments.}}}{
%	We are grateful to the referees for their stimulating comments and suggestions. We acknowledge partial support by the %FRA2018 grant `Uncertainty Modelling: Mathematical Methods and Applications'.
%}
%
%\bibliographystyle{model1-num-names}
%\bibliographystyle{plain}
%\section*{\refname}
%\printbibliography
%\small

\bibliographystyle{plain}
\bibliography{NL_Dilation_Arxiv}{}

\end{document}